\documentclass[11pt]{article}

\usepackage[utf8]{inputenc}
\usepackage[english]{babel}
\usepackage{amsmath}
\usepackage{amsthm}
\usepackage{amssymb}
\usepackage{blkarray}
\usepackage{appendix}

\newtheorem{theorem}{Theorem}
\newtheorem{corollary}{Corollary}
\newtheorem{lemma}{Lemma}
\newtheorem{definition}{Definition}

\newtheorem{example}{Example}
\newtheorem{conjecture}{Conjecture}
\newcommand{\oN}{\mathbb N}
\newcommand{\oR}{\mathbb R}
\newcommand{\ignore}[1]{}
\newtheorem{question}{Question}

\date{\today}

\usepackage[utf8]{inputenc}
\usepackage[utf8]{inputenc}
\usepackage[english]{babel}
\usepackage{amsthm}
\usepackage{amsmath}
\usepackage{color}
\usepackage[english]{babel}
\usepackage[hidelinks]{hyperref}
\usepackage{amsmath,amssymb,amsthm}
\usepackage{mathtools}
\usepackage{graphicx}
\usepackage{dsfont}

\usepackage{float}

\newcommand{\EnL}{E}
\newcommand{\FmL}{F}

\newcommand{\uma}{\min\{i,j\}}

\newcommand{\comba}{a(}
\newcommand{\combb}{b(u,}
\newcommand{\combc}{c(u)}

\newcommand{\g}{g(u,z)}

\newcommand{\MP}{{\mathcal P}}

\newcommand{\Sym}{\text{\rm Sym}}

\usepackage{pgf,tikz,pgfplots}
\usepackage{mathrsfs}
\usetikzlibrary{arrows}
	\tikzset{
	triple/.style args={[#1] in [#2] in [#3]}{
		#1,preaction={preaction={draw,#3},draw,#2}
}}

\definecolor{xdxdff}{rgb}{0.49019607843137253,0.49019607843137253,1.}
\definecolor{ududff}{rgb}{0.30196078431372547,0.30196078431372547,1.}

\title{Optimizing  hypergraph-based  polynomials  modeling job-occupancy in queueing with redundancy scheduling}

\author{
		Daniel Brosch \thanks{Tilburg University, \url{d.brosch@uvt.nl}}		
        		\and
		Monique Laurent \thanks{Centrum Wiskunde \& Informatica (CWI), Amsterdam, and Tilburg University,
   	\url{Monique.Laurent@cwi.nl}}		
			\and
        Andries Steenkamp \thanks{Centrum Wiskunde \& Informatica (CWI), Amsterdam,
        \url{Andries.Steenkamp@cwi.nl}}
     }

	\date{\today}

\begin{document}
\maketitle

\begin{abstract}
We investigate  two classes of multivariate polynomials with  variables indexed by the edges of a  uniform hypergraph and coefficients depending on certain patterns of union of edges. These polynomials arise naturally to model job-occupancy in some queuing problems 
with redundancy scheduling policy.  The question, posed by Cardinaels, Borst and van Leeuwaarden (arXiv:2005.14566, 2020), is to decide whether their global minimum over the standard simplex is attained at the uniform probability distribution.
By exploiting symmetry properties of these polynomials we can give a positive answer for the first class and partial results for the second one, where we in fact show a stronger convexity property of these polynomials over the simplex.

\end{abstract}

\section{Introduction}\label{secsetting}

We consider the minimization of a class of polynomials over the standard simplex. These polynomials have their variables labelled by the  edges of a complete uniform hypergraph and their coefficients are defined in terms of some cardinality patterns of unions of edges. They arise naturally within the modelling of job-occupancy in some queuing problems with redundancy scheduling policies 
\cite{CBS}. The question is whether these polynomials attain their minimum value at the barycenter of the standard simplex, which corresponds to showing optimality of the uniform distribution for the underlying queuing problem.
This paper is devoted to this question.

We now  introduce the classes of polynomials of interest.
Given integers $n,L\ge 2$ we set $V=[n]=\{1,\ldots,n\}$ and
$E=\{e\subseteq V: |e|=L\}$, so that $(V,E)$ can be seen as the complete $L$-uniform hypergraph on $n$ elements.
We set $m:=|E|={n\choose L}$, where we omit the explicit dependence on $n,L$ to simplify notation,
and we let $$\Delta_m=\Big\{x=(x_e)_{e\in E}\in \oR^m: x\ge 0, \sum_{e\in E}x_e=1\Big\}$$ denote the  standard simplex in $\oR^m$. The elements of $\Delta_m$ correspond to probability vectors on $m$ items and  the barycenter $x^*={1\over m}(1,\ldots,1)$ of $\Delta_m$ corresponds to the uniform probability vector.

Given an integer $d\ge 2$ we consider the following $m$-variate polynomial in the variables $x=(x_e: e\in E)$, which is a main player in the paper:
\begin{equation}\label{eqpd0}
p_d(x)=\sum_{(e_1,\ldots, e_d)\in E^d} {1\over |e_1\cup \ldots \cup e_d|} x_{e_1}\cdots x_{e_d}.
\end{equation}
So $p_d$ is homogeneous with degree $d$.
We are interested in  the following optimization problem
\begin{equation*}\label{eqopt}
p_d^*:=\min_{x\in \Delta_m} p_d(x),
\end{equation*}
asking to minimize the polynomial $p_d$ over the simplex $\Delta_m$.
Our main objective is to show that the global minimum is attained at  the uniform probability vector $x^*$.
The following is the main result of the paper.

\begin{theorem}\label{theomainpd}
The global minimum of the polynomial $p_d$ from (\ref{eqpd0}) over the standard simplex $\Delta_m$ is attained at the barycenter $x^*={1\over m}(1,\ldots,1)$ of $\Delta_m$.
\end{theorem}

There is a second related class of polynomials of interest
\begin{equation}\label{eqfd0}
f_d(x) = \sum_{(e_1,\ldots,e_d)\in E^d} \prod_{i=1}^d {x_{e_i}\over |e_1\cup \ldots \cup e_i|},
\end{equation}
 also homogeneous with degree $d$.
Note that for  $d=1$ both polynomials coincide: $p_1=f_1$ and, for $d=2$, we have $f_2= {1\over L} p_2$. Here too the question is whether  the minimum of $f_d$ over the standard simplex $\Delta_m$  is attained at the uniform probability vector $x^*$.

\begin{question}\label{question}
Given integers $n,d,L\ge 2$ is it true that the polynomial $f_d(x)$ in (\ref{eqfd0}) attains its minimum over $\Delta_m$ at the barycenter $x^*$ of $\Delta_m$?
\end{question}
 As noted above the answer is positive for $d=2$ and $n,L\ge 2$. As a further partial result we give
 a positive answer for the case of degree $d=3$ and edge size $L=2$.

 \begin{theorem}\label{theomainfd3}
For  $d= 3$ and  $L=2$ the global minimum of the polynomial $f_d$ from (\ref{eqfd0}) over the standard simplex $\Delta_m$ is attained at the barycenter $x^*={1\over m}(1,\ldots,1)$ of $\Delta_m$.
\end{theorem}

As we will mention in the next section the above question about the polynomials $f_d$,  posed by the authors of \cite{CBS} (in the case $L=2$), is motivated by its relevance to a problem in queueing theory. The polynomials $f_d$ can be seen as a variant of the  polynomials $p_d$ and as mentioned above  both classes   coincide (up to scaling) for degree $d=2$. For the polynomials $p_d$ we can give a full answer and show that they indeed attain their minimum at the barycenter of $\Delta_m$. The analysis of the polynomials $f_d$ is  technically much more involved and we have only partial results so far. In both cases  the key ingredient is showing that the polynomials are convex, i.e., that they have positive semidefinite Hessians. 
 It turns out that the Hessian of the polynomial $p_d$ enters in some  way as a component of  the Hessian of the polynomial $f_d$. So this forms a natural motivation for the study of the polynomials $p_d$, though they form a natural class of symmetric polynomials that are interesting for their own sake.

Exploiting symmetry plays a central role in our proofs. Indeed the key idea is to show that the polynomials are convex, which, combined with their symmetry properties, implies that the global minimum is attained at the barycenter of the simplex. For this we show that their Hessian matrices are positive semidefinite at each point of the simplex, which we do  through  exploiting again  their symmetry structure and links to Terwilliger algebras.

Symmetry is a widely used ingredient in optimization, in particular  in semidefinite optimization and algebraic questions involving polynomials. We mention a few  landmark examples  as background information. Symmetry can indeed be used to formulate equivalent, more compact reformulations for semidefinite programs. The underlying mathematical fact is Wedderburn-Artin theory, which shows that  matrix $*$-algebras can be block-diagonalized (see Theorem \ref{theartwed} {below). 
An early well-known example is the linear programming reformulation from \cite{Schrijver1979} for the Lov\'asz theta number 
of Hamming graphs, showing the link to the Delsarte bound and Bose-Mesner algebras of Hamming schemes \cite{Delsarte73,DelsarteLevenshtein}.
Symmetry is  used more generally to give tractable reformulations for the semidefinite bounds arising from the next levels of Lasserre's hierarchy in  \cite{Schrijver2005} (which gives the explicit block-diagonalization for the Terwilliger algebra of Hamming schemes, see Theorem \ref{schrter} below) and, e.g.,  in  \cite{GST}, \cite{GMS12},   \cite{Laurent2007}, \cite{LPS2017}. For more examples and a broad exposition  about the use of symmetry in semidefinite programming
we refer, e.g., to \cite{BGSV,dKPS}  and further references therein.  Symmetry is also a crucial ingredient in the study of algebraic questions about polynomials, like representations in terms of  sums of squares, and in polynomial optimization.
 We refer to \cite{GatermanParrilo} for a broad exposition and, e.g.,  to  \cite{RTAL} (for compact reformulations of Lasserre relaxations of symmetric polynomial optimization problems), \cite{Riener} 
(for methods to reduce the number of variables in programs involving symmetric polynomials), and the recent works \cite{RSST2018,RSTTuran2018} (which consider symmetric polynomials with variables indexed by the $k$-subsets hypercube (as in our case) and uncover links with the theory of flag algebras by Razborov \cite{Razborov}).



\begin{example}
As an illustration let us consider  the polynomial $p_d$ for edge size  $L=2$. 
For $d=1$, we have $p_1(x)= {1\over 2}\sum_{e\in E}x_e$. For $d=2$ we have
$$p_2(x)= {1\over 2} \sum_{e\in E} x_e^2+ {1\over 3}\sum_{{(e_1,e_2)\in E^2:}\atop{ |e_1\cup e_2|=3}}x_{e_1}x_{e_2}+{1\over 4}
\sum_{{(e_1,e_2)\in E^2:}\atop{|e_1\cup e_2|=4}}x_{e_1}x_{e_2}.$$
Using the notation $p_d(x)=\sum_{\underline{e}=(e_1,\ldots,e_d)\in E^d} c_{\underline{e}}x_{e_1}\cdots x_{e_d}$ from relation (\ref{eqpd}) below for the polynomial $p_d$, we show  in Figure \ref{figEdgesTup2}
the three possible patterns for pairs of edges $\underline{e}=(e_1,e_2)$ and the corresponding coefficients $c_{\underline{e}}$.

\begin{figure}[H]
	\centering
	\begin{tikzpicture}[line cap=round,line join=round,>=triangle 45,x=1.7cm,y=1.7cm]
	\draw [double,line width=2.pt] (-3.,0.5)-- (-2.,-0.5);
	
	\draw [line width=2.pt] (-1.,0.5)-- (0.,-0.5);
	\draw [line width=2.pt] (0.,-0.5)-- (1.,0.5);
	
	\draw [line width=2.pt] (2.,0.5)-- (3.,-0.5);
	\draw [line width=2.pt] (3.,0.5)-- (4.,-0.5);
	
	\begin{scriptsize}
	\draw [fill=ududff] (-3.,0.5) circle (4pt);
	\draw [fill=ududff] (-2.,-0.5) circle (4pt);

	\draw [fill=ududff] (-1.,0.5) circle (4pt);
	\draw [fill=xdxdff] (0.,-0.5) circle (4pt);
	\draw [fill=ududff] (1.,0.5) circle (4pt);

	\draw [fill=ududff] (2.,0.5) circle (4pt);
	\draw [fill=ududff] (3.,-0.5) circle (4pt);
	\draw [fill=ududff] (3.,0.5) circle (4pt);
	\draw [fill=ududff] (4.,-0.5) circle (4pt);

	\draw[color=black] (-2.9,-.5) node {$c_{\underline{e}} = \frac{1}{2}$};
	\draw[color=black] (-.9,-.5) node {$c_{\underline{e}} = \frac{1}{3}$};
	\draw[color=black] (2.1,-.5) node {$c_{\underline{e}} = \frac{1}{4}$};
	
	\end{scriptsize}
	\end{tikzpicture}
	\caption{The three patterns of pairs of edges in case ($d=2, L=2$)}
	\label{figEdgesTup2}
\end{figure}
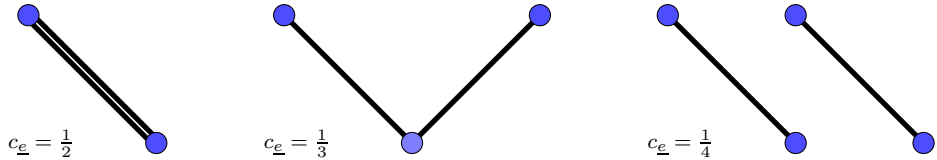

 In the same way, for $d\ge 3$,
$p_d(x)=\sum_{k=2}^{2d} {1\over k} q_{d,k}(x)$, where
the summand $q_{d,k}(x)$
is a summation over all $d$-tuples of edges with a given pattern, depending on the cardinality  of their union:
$$q_{d,k}(x)= \sum_{{(e_1,\ldots,e_d)\in E^d:}\atop{|e_1\cup \ldots \cup e_d|=k}}x_{e_1}\cdots x_{e_d}.$$
For the case $d=3$ we need to consider the values $k=2,3,4,5,6$; as an illustration we show in Figure \ref{figEdgesTup3} all the possible patterns of triplets of edges $\underline{e}=(e_1,e_2,e_3)$ and the corresponding coefficients $c_{\underline{e}}$ that contribute to the summands $q_{3,k}$.
\end{example}

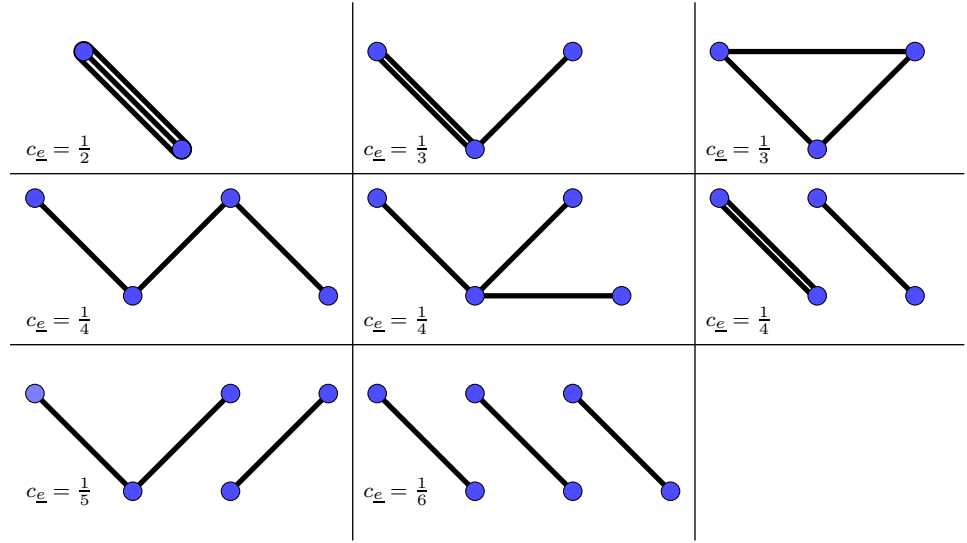
\begin{figure}[H]
	\centering
	\begin{tikzpicture}[line cap=round,line join=round,>=triangle 45,x=0.65cm,y=0.65cm]
	\draw [triple = {[line width=2.pt] in
		[line width=4.pt,white] in
		[line width=8.pt]}] (-6.,6.)-- (-4.,4.);
	
	\draw [double,line width=2.pt] (0.,6.)-- (2.,4.);
	\draw [line width=2.pt] (2.,4.)-- (4.,6.);
	
	\draw [line width=2.pt] (7.,6.)-- (9.,4.);
	\draw [line width=2.pt] (9.,4.)-- (11.,6.);
	\draw [line width=2.pt] (7.,6.)-- (11.,6.);
	
	\draw [line width=2.pt] (-7.,3.)-- (-5.,1.);
	\draw [line width=2.pt] (-5.,1.)-- (-3.,3.);
	\draw [line width=2.pt] (-3.,3.)-- (-1.,1.);
	\draw [line width=2.pt] (5.,1.)-- (2.,1.);
	\draw [line width=2.pt] (0.,3.)-- (2.,1.);		
	\draw [line width=2.pt] (4.,3.)-- (2.,1.);
	\draw [double,line width=2.pt] (7.,3.)-- (9.,1.);
	\draw [line width=2.pt] (9.,3.)-- (11.,1.);
	
	\draw [line width=2.pt] (-7.,-1.)-- (-5.,-3.);
	\draw [line width=2.pt] (-5.,-3.)-- (-3.,-1.);
	\draw [line width=2.pt] (-3.,-3.)-- (-1.,-1.);
	
	\draw [line width=2.pt] (0.,-1.)-- (2.,-3.);
	\draw [line width=2.pt] (2.,-1.)-- (4.,-3.);
	\draw [line width=2.pt] (4.,-1.)-- (6.,-3.);

	\draw [line width= 0.5pt] (-0.5,7.)-- (-0.5,-4.);
	\draw [line width= 0.5pt] (6.5,7.)-- (6.5,-4.);
	
	\draw [line width= 0.5pt] (-7.5,3.5)-- (12,3.5);
	\draw [line width= 0.5pt] (-7.5,0)-- (12,0);

	\begin{scriptsize}
	\draw [fill=ududff] (-6.,6.) circle (3.5pt);
	\draw [fill=ududff] (-4.,4.) circle (3.5pt);
	\draw [fill=ududff] (0.,6.) circle (3.5pt);
	\draw [fill=ududff] (2.,4.) circle (3.5pt);
	\draw [fill=ududff] (4.,6.) circle (3.5pt);

	\draw [fill=ududff] (7.,6.) circle (3.5pt);
	\draw [fill=ududff] (9.,4.) circle (3.5pt);
	\draw [fill=ududff] (11.,6.) circle (3.5pt);
	
	\draw [fill=ududff] (-7.,3.) circle (3.5pt);
	\draw [fill=ududff] (-5.,1.) circle (3.5pt);
	\draw [fill=ududff] (-3.,3.) circle (3.5pt);
	\draw [fill=ududff] (-1.,1.) circle (3.5pt);

	\draw [fill=ududff] (0.,3.) circle (3.5pt);
	\draw [fill=ududff] (4.,3.) circle (3.5pt);
	\draw [fill=ududff] (2.,1.) circle (3.5pt);
	\draw [fill=ududff] (5.,1.) circle (3.5pt);

	\draw [fill=ududff] (7.,3.) circle (3.5pt);
	\draw [fill=ududff] (9.,1.) circle (3.5pt);
	\draw [fill=ududff] (9.,3.) circle (3.5pt);
	\draw [fill=ududff] (11.,1.) circle (3.5pt);
	
	\draw [fill=xdxdff] (-7.,-1) circle (3.5pt);	
	\draw [fill=ududff] (-5.,-3) circle (3.5pt);
	\draw [fill=ududff] (-3.,-1) circle (3.5pt);
	\draw [fill=ududff] (-3.,-3) circle (3.5pt);
	\draw [fill=ududff] (-1.,-1) circle (3.5pt);
	
	\draw [fill=ududff] (0.,-1.) circle (3.5pt);
	\draw [fill=ududff] (2.,-3.) circle (3.5pt);
	\draw [fill=ududff] (2.,-1.) circle (3.5pt);
	\draw [fill=ududff] (4.,-3.) circle (3.5pt);
	\draw [fill=ududff] (4.,-1.) circle (3.5pt);
	\draw [fill=ududff] (6.,-3.) circle (3.5pt);

	\draw[color=black] (-6.5,4.0) node {$c_{\underline{e}} = \frac{1}{2}$};
	\draw[color=black] (0.4,4.0) node {$c_{\underline{e}} = \frac{1}{3}$};
	\draw[color=black] (7.4,4.0) node {$c_{\underline{e}} = \frac{1}{3}$};  
	
	\draw[color=black] (-6.5,0.5) node {$c_{\underline{e}} = \frac{1}{4}$};
	\draw[color=black] (0.4,0.5) node {$c_{\underline{e}} = \frac{1}{4}$};
	\draw[color=black] (7.4,0.5) node {$c_{\underline{e}} = \frac{1}{4}$};
	
	\draw[color=black] (-6.5,-3.0) node {$c_{\underline{e}} = \frac{1}{5} $};
	\draw[color=black] (0.4,-3.0) node {$c_{\underline{e}} = \frac{1}{6}$};
	\end{scriptsize}
	
	\end{tikzpicture}
	\caption{The eight patterns of triplets of edges in case ($d=3, L=2$)}
	\label{figEdgesTup3}
\end{figure}

\medskip\noindent
{\bf Organization of the paper.} 
In the rest of this section  we first indicate in Section \ref{secmotivation} how the polynomials $f_d$  naturally arise within a problem of  queuing theory with redundancy scheduling policies.  After that we present in Section \ref{secsketch} the main ideas of the proofs,
which highly rely on exploiting symmetry properties of the polynomials.
 This involves in particular using the Terwilliger algebra of the binary Hamming cube, so we  include  some preliminaries about  these Terwilliger algebras in Section \ref{secPrelimTer}.

In Section \ref{secprooftheopd} we give the full proof for Theorem~\ref{theomainpd} showing that the polynomials $p_d$ attain their global minimum at the barycenter of the simplex and, in Section \ref{secprooffd}, we investigate the second class of polynomials $f_d$.
We prove several properties of these polynomials, which we use to show Theorem~\ref{theomainfd3}. 
We also present a range of values of $(n,d,L)$ for which the polynomials $f_d$ are indeed convex and thus Question \ref{question} has a positive answer.

\medskip\noindent
{\bf Some notation.}
Throughout we let $I,J$ denote the identity matrix and the all-ones matrix, whose size should be clear from the context. When we want to specify the size we let $I_n$ (resp., $J_n$) denote the $n\times n $ identity matrix (resp., all-ones matrix) and, given two integers $n,m\ge 1$, $J_{m,n}$ denotes the $m\times n$ all-ones matrix. For a symmetric matrix $A$ the notation $A\succeq 0$ means that $A$ is positive semidefinite. 
Given two matrices $A,B \in \oR^{n\times n}$ we let $A\circ B\in \oR^{n\times n}$ denote their Hadamard product, with entries
$(A\circ B)_{ij} =A_{ij}B_{ij}$ for $i,j\in [n]$. It is known that $A\succeq 0$ and $B\succeq 0$ implies $A\circ B\succeq 0$.

For a sequence $\alpha\in\oN^n$ we set $|\alpha|=\sum_{i=1}^n \alpha_i$ and, for an integer $d\in \oN$, we set $\oN^n_d=\{\alpha\in\oN^n: |\alpha|=d\}$.
Given a vector $x\in\oR^n$ and $\alpha\in \oN^n$ we set $x^\alpha=x_1^{\alpha_1}\cdots x_n^{\alpha_n}$.
Throughout we let $u_1,\ldots,u_m$ denote the standard basis of $\oR^m$, where all entries of $u_i$ are 0 except its $i$-th entry  which is equal to 1.
We let $\Sym(n)$ denote the set of permutations of the set $V=[n]$.

\subsection{Motivation }\label{secmotivation}

Our motivation for the study of the  polynomials $p_d$ and $f_d$ comes from their relevance to a problem in queueing theory.
The question whether they attain their minimum at the uniform probability distribution was posed to us by the authors of \cite{CBS}, who use a positive answer to this question to establish a result about the asymptotic behaviour of the job occupancy in a parallel-servers system with redundancy scheduling in the light-traffic regime.
In what follows we will give only a high level sketch of this connection and we refer to the paper \cite{CBS} for a detailed exposition.
We also refer to \cite{CBS} for an extended review of the relevant literature.

A crucial mechanism that has been considered  to improve the performance of parallel-servers systems in queueing theory is {redundancy scheduling}.  The key feature of this policy is that several  replicas  are created for each arriving job, which are then assigned to distinct servers (and then, as soon as the first of these replicas completes (or enters) service on a server the remaining ones are stopped).
The underlying idea is that sending replicas of the same job to several servers will  increase the chance of having shorter queueing times. This however must be weighted against the risk of wastage of capacity. An important question is thus to assess the impact of redundancy scheduling policies. While most papers in the literature of redundant scheduling assume that the set of servers to which the replicas are sent is selected uniformly at random, the paper \cite{CBS} considers the case when the set of servers is selected according to a given probability distribution and  it investigates what is the impact of this probability distribution on the performance of the system.
It is shown there that while  the impact remains relatively limited in the heavy-traffic regime,  the system occupancy is much more sensitive to the selected probability distribution in the light-traffic regime.



\medskip
We will now only introduce a few elements of the model considered in \cite{CBS}, so that we can make the link to the polynomials studied in this paper. We keep our presentation high level and refer to \cite{CBS} for details. The setting is as follows. There are $n$ parallel servers,  with average speed $\mu$.
Jobs arrive as a Poisson process of rate $n\lambda$ for some $\lambda > 0$.
When a job arrives, $L$ replicas of it are created that are sent - with probability $x_e$ - to a subset $e\subseteq [n]$ of $L$ servers.
Here, $L\ge 2$ is an integer and $x=(x_e)_{e\in E}$ is a probability distribution  on the set $E=\{e\subseteq [n]: |e|=L\}$ of possible collections of $L$ servers.
As noted in \cite{CBS}  this can be seen as selecting  an edge $e\in E$ with probability $x_e$ in the uniform hypergraph $(V=[n],E)$ (with edge size $L$).

An important performance parameter is the system occupancy at time $t$, which is represented by a vector $(e_1,..., e_{M}) \in {E}^{M}$,
where $M=M(t)$ is the total number of jobs present in the system and $e_i \in {E}$ is the collection of servers  to which the replicas of the $i$-th longest job in the system have been assigned. 
Under suitable stability conditions and assuming each server has the same speed $\mu$ and service requirements of the jobs are independent and exponentially distributed with unit mean,
the stationary distribution of the  occupancy of this edge selection is given by
$$\pi(e_1,\ldots,e_M)=C \prod_{i=1}^M {n\lambda x_{e_i}\over \mu |e_1\cup \ldots \cup e_i|}$$
for some constant $C>0$ (\cite{GZDHBHSW}, see relation (3) in \cite{CBS}).
Following \cite{CBS}, let  $Q_\lambda(x)$ be a random variable with the stationary distribution of the system occupancy when the edge selection is given by the probability vector $x=(x_e)_{e\in E}$. It then follows that, for any integer $d\ge 1$, the probability that $d$ jobs are present in the system is given by
$$\mathbb {P}\{Q_\lambda(x)= d\}=\sum_{(e_1,\ldots,e_d)\in E^d} \pi(e_1,\ldots,e_d).$$
Hence, $\mathbb {P}\{Q_\lambda(x)= 0\}=C$ and
$$\mathbb {P}\{Q_\lambda(x)= d\}=\mathbb {P}\{Q_\lambda(x)= 0\} \Big({n\lambda\over \mu}\Big)^d \sum_{(e_1,\ldots,e_d)\in E^d}
\prod_{i=1}^d {x_{e_i}\over |e_1\cup \ldots \cup e_i|}.
$$
(See relation (11) in \cite{CBS}). Therefore, $\mathbb {P}\{Q_\lambda(x)= d\}$ is the polynomial $f_d(x)$ (up to a scalar multiple).
In \cite{CBS}  the light-traffic regime is considered, i.e., when $\lambda \downarrow 0$, in the case $L=2$.  By doing a Taylor expansion one can see that
$$\mathbb {P}\{Q_\lambda(x)= 0\}= 1+o(1),\
\mathbb {P}\{Q_\lambda(x)\ge  d\}= \Big({n\lambda\over \mu}\Big)^d f_d(x)+o(\lambda^d)$$
(see relation (13) in \cite{CBS}).
Therefore, with $x^*=(1,\ldots,1)/|E|$ denoting the uniform probability vector, we have
$$\lim_{\lambda\downarrow 0} {\mathbb {P}\{Q_\lambda(x^*)\ge  d\}\over \mathbb {P}\{Q_\lambda(x)\ge  d\}}
= \lim_{\lambda\downarrow 0} {f_d(x^*)+o(1)\over f_d(x)+o(1)}.$$
Hence, if the polynomial $f_d$ attains its minimum at the uniform distribution  $x^*$, then one has
$$\lim_{\lambda\downarrow 0} {\mathbb {P}\{Q_\lambda(x^*)\ge  d\}\over \mathbb {P}\{Q_\lambda(x)\ge  d\}} \le 1.$$
This indicates that in the light-traffic regime the system occupancy is minimized when selecting uniformely at random the assignments to the servers of the job replicas. This thus motivates Question \ref{question} of showing that the polynomial $f_d$ attains its minimum over the probability simplex at the uniform point $x^*$.

\subsection{Sketch of proof}\label{secsketch}


Here we give a sketch of proof for our main results. We start with indicating the main steps for proving Theorem \ref{theomainpd}, dealing with the class of polynomials $p_d$ and after that we briefly indicate how to deal with the polynomials $f_d$.

A first easy observation is that in order to show that the polynomial $p_d$ attains its minimum at the barycenter of the standard simplex $\Delta_m$  it suffices to show that  $p_d$ is convex over $\Delta_m$. This follows from a symmetry argument, namely we exploit the fact   that the polynomial $p_d$ is invariant under the permutations of the edge set $E$ that are induced by permutations of $[n]$.

\begin{lemma}\label{lemcondition}
Assume  the polynomial $p_d$ is convex on the simplex $\Delta_m$. Then the point $x^*=(1/m)(1,\ldots,1)\in \Delta_m$ is a global minimizer of $p_d$ over $\Delta_m$.
\end{lemma}

\begin{proof} 
The key  fact we use is that the polynomial $p_d$ enjoys some symmetry property; namely, for any tuple $(e_1,\ldots,e_d)\in E^d$,   the coefficient of the monomial $x_{e_{1}}\cdots x_{e_{d}}$ in $p_d$ is $1/|e_{1}\cup\ldots \cup e_{d}|$, which depends only on the cardinality of the set $e_{1}\cup\ldots\cup e_{d}$.
Recall that $E=\{e\subseteq V=[n]: |e|=L\}$. Any permutation $\sigma\in \Sym(n)$ of $[n]$ induces a permutation of $E$ (still denoted $\sigma$)  by setting $\sigma(e)=\{j_{\sigma(1)},\ldots,j_{\sigma(L)}\}$ for  $e=\{j_1,\ldots,j_L\}\in E$. In turn, $\sigma$ acts on $\Delta_m$ by setting
$\sigma(x)=(x_{\sigma(e)})_{e\in \Delta_m}$ for $x=(x_e)_{e\in E}\in \Delta_m$.
We now observe that $p_d$ is invariant under this action of  permutations $\sigma\in \Sym(n)$. Indeed, for any $\sigma\in \Sym(n)$,  we have
$$\begin{array}{ll}
\sigma(p_d)(x)=p_d(\sigma(x))& =\sum_{(e_1,\ldots,e_d)\in E^d}{1\over |e_1\cup \ldots \cup e_d|} x_{\sigma(e_1)}\cdots x_{\sigma(e_d)}\\
&= \sum_{(f_1,\ldots,f_d)\in E^d} {1\over |\sigma^{-1}(f_1)\cup \ldots \cup \sigma^{-1}(f_d)|} x_{f_1}\cdots x_{f_d}\\
&= \sum_{(f_1,\ldots,f_d)\in E^d} {1\over |f_1\cup \ldots \cup f_d|}x_{f_1}\cdots x_{f_d} \\
&= p_d(x).
\end{array}$$
Let $x\in \Delta_m$ be a global minimizer of $p_d$. For any permutation $\sigma \in \Sym(n)$ the permuted point $\sigma(x)$ belongs to $\Delta_m$ and $p_d(x)=p_d(\sigma(x))$ holds.
Hence, for  the full symmetrization of $x$,
$$x^*:={1\over n!}\sum_{\sigma \in\Sym(n)} \sigma(x),$$
we have $x^*\in \Delta_m$ and
$$p_d(x^*)\le {1\over n!}\sum_{\sigma\in \Sym(n)}p_d(\sigma(x))= p_d(x),$$
where the inequality holds since $p_d$ is  convex over $\Delta_m$.
This shows that $x^*$ is again a global minimizer of $p_d$ in $\Delta_m$.
It suffices now to observe that, by construction, $x^* = (1/m)(1,\ldots,1)$.
\end{proof}


Therefore we are left with the task of showing that the polynomial $p_d$ is convex over the simplex $\Delta_m$ or,
equivalently,  that its Hessan matrix
$$H(p_d)(x)=(\partial^2 p_d(x)/\partial x_e \partial x_f)_{e,f\in E}$$
 is positive semidefinite over $\Delta_m$. This forms the core technical part of the proof. Here is a rough sketch of our proof technique.

A first step is to express the Hessian matrix  as a matrix polynomial, involving a collection of matrices $M_\gamma$;
see Lemma \ref{lemHpd}.
The next step is to show that each of the matrices $M_\gamma$ appearing in this decomposition of the Hessian is positive semidefinite.
For this, we reduce to the task of showing that a certain set of well-structured matrices are positive semidefinite, see  Lemmas \ref{proppdconvex} and \ref{symmRedu}.
This last task is done by showing that these matrices lie in the Terwilliger algebra of the Hamming cube, which enables us to exploit its explicitly known block-diagonalization.  The proof is then concluded by using an integral representation argument, see Section \ref{secdge3Lge3}.

\medskip
The treatment for the polynomials $f_d$ has the same starting point:  the polynomial $f_d$ is invariant under any permutation of the edge set $E$ induced by permutations of $[n]$  and thus it suffices to show that $f_d$ is convex in order to conclude that it attains its global minimum at the barycenter of the simplex (i.e., the analogue of Lemma \ref{lemcondition} holds for $f_d$).
After that we again express the Hessian matrix $H(f_d)$ as a matrix polynomial, involving a collection of marices $Q_\gamma$; see Lemma \ref{lemHfd}. Hence here too the  task boils down to showing that each of these matrices $Q_\gamma$  is positive semidefinite.
This task turns out to be considerably more difficult than for the matrices $M_\gamma$ which occurred in the analysis of the polynomial $p_d$.
As a first step toward the analysis of the matrices $Q_\gamma$ we give a recursive reformulation for them, which also makes apparent how the matrices $M_\gamma$  enter their definition (namely as a factor of a Hadamard product definition of $Q_\gamma$); see Lemma~\ref{lembahat}.
Based on this we can show that the matrices $Q_\gamma$ are indeed positive semidefinite in the case $d=3$ and $L=2$, thus showing Theorem \ref{theomainfd3}; see Section \ref{secfd3}.

\subsection{Preliminaries on the Terwilliger algebra}\label{secPrelimTer}

As mentioned above we need to exploit the symmetry structure of the polynomial $p_d$ in order to show that its Hessian matrix is positive semidefinite.
A crucial ingredient will be that the Hessiam matrix can be decomposed into matrices that (after some reduction steps) all lie in the Terwiliger algebra of the binary Hamming cube.
We begin with introducing the definition of the Terwiliger algebra $\mathcal A_n$ of the binary Hamming cube on $n$ elements.


\begin{definition}[\textbf{Terwilliger algebra of the binary Hamming cube}]
Let $\MP_n$ denote  the collection of all subsets of the set $V = [n]$. For every triple of nonnegative integers $i,j,t$ we define the $2^n\times 2^n$ matrix $D_{i,j}^t$, indexed by $\MP_n$, with entries
$$
\left(D_{i,j}^t\right)_{S,T} =
\begin{cases}
      1 & if ~ |S| = i, |T|= j, |S\cap T| = t,\\
      0 & \text{ else}
\end{cases}.
$$
for sets $S,T\in\MP_n$.
Then the {\em  Terwilliger algebra of the binary Hamming cube}, denoted by  $\mathcal A_n$, is defined as the (real) span of all these matrices:
$$\mathcal A_n=\Big\{\sum_{i,j,t\ge 0} x^t_{i,j}D^t_{i,j}: x^t_{i,j}\in \oR\Big\}.$$
\end{definition}

It is easy to see that $\mathcal A_n$ is a \emph{matrix $\ast$-algebra}, i.e., $\mathcal A_n$ is  closed under taking  linear combinations, matrix multiplications and transposition. One way to see this is by realizing that the matrices $D_{i,j}^t$ are exactly the indicator matrices of the orbits of pairs in $\MP_n\times \MP_n$ under the element-wise action of the symmetric group $\Sym(n)$.

%

All matrix $\ast$-algebras can be block-diagonalized by Artin-Wedderburn theory (see \cite{ArtinWedderburn}, see also \cite{BEK} for a proof).

\begin{theorem}[Artin-Wedderburn]\label{theartwed}
Let $\mathcal A$ be a matrix $\ast$-algebra. Then there exist nonnegative integers $d$ and $m_1,\ldots,m_d$ and  a $\ast$-algebra isomorphism
$$\varphi\colon \mathcal A\to\bigoplus_{k=1}^d \mathbb{C}^{m_k\times m_k}.$$
\end{theorem}
The important property here   is that $\varphi$ is an algebra isomorphism. Hence we know that this isomorphism maintains positive semidefiniteness: for any matrix $A\in \mathcal A$, we have
$A\succeq 0 \Longleftrightarrow \varphi (A)\succeq 0$.
Moreover, 
the matrix $\varphi(A)$ is block-diagonal, with $d$ diagnal blocks of sizes $m_1,\ldots,m_d$.  This is a crucial property which can be exploited in order
to get a more efficient way of encoding positive semidefiniteness of matrices in $\mathcal A$.

The explicit block-diagonalization of the Terwilliger algebra $\mathcal A_n$ was given  by Schrijver \cite{Schrijver2005}.

\begin{theorem}[Schrijver \cite{Schrijver2005}]\label{schrter}
The Terwiliger algebra   $\mathcal A_n$ can be block-diagonalized into $\lfloor \frac{n}{2}\rfloor +1$ blocks, of sizes $m_k = n-2k+1$ for $k=0,\ldots,\lfloor \frac{n}{2}\rfloor$. The algebra isomorphism $\varphi$ sends the matrix
  $$A=\sum_{i,j,t=0}^{n}x_{i,j}^tD_{i,j}^t$$
  to the block-matrix
  $\varphi(A)=\oplus_{k=0}^{\lceil n/2\rceil} B_k$, where the matrix  $B_k\in\oR^{m_k\times m_k}$ is given by
  \begin{equation}\label{eqBk}
  B_k := \left(\binom{n-2k}{i-k}^{-\frac{1}{2}}\binom{n-2k}{j-k}^{-\frac{1}{2}}\sum_{t}\beta_{i,j,k}^t x_{i,j}^t\right)_{i,j=k}^{n-k}
  \end{equation}
  for $k=0,1,\ldots,\lfloor \frac{n}{2}\rfloor.$
Here, for any nonnegative integers $i,j,t,k$,  we set
  \begin{equation}\label{eqbeta}
  \beta_{i,j,k}^t := \sum_{u=0}^{n}(-1)^{u-t}\binom{u}{t}\binom{n-2k}{n-k-u}\binom{n-k-u}{i-u}\binom{n-k-u}{j-u}.
  \end{equation}
  In particular we have
  \begin{equation}\label{eqequiv}
  \sum_{i,j,t=0}^{n}x_{i,j}^tD_{i,j}^t\succeq 0 \Longleftrightarrow  B_k\succeq 0\ \text{ for }  k=0,1,\ldots,\lfloor \frac{n}{2}\rfloor.
  \end{equation}
  \end{theorem}

\section{Proof of Theorem \ref{theomainpd}}\label{secprooftheopd}
In this section we give the proof of  Theorem \ref{theomainpd}.
As a warm-up we start with the special case when the degree is $d=2$ and the edge size is $L=2$, where we can easily show that the polynomial $p_2$ is convex.

After that we proceed to the general case.
We follow the steps as sketched above: first we express the Hessian matrix of $p_d$ as a matrix polynomial and we indicate some reductions that lead to the task of showing that a set of well-structured matrices are positive semidefinite. After that we show the positive semidefiniteness of these matrices by exploiting a link to the Terwilliger algebra of the Boolean Hamming cube.

\subsection{The case $d=2$ and $L=2$}\label{secd2L2}

Here  we consider the polynomial
$$p_2(x)=\sum_{e,f\in E} {1\over |e\cup f|} x_ex_f,$$
where $E=\{e\subseteq [n]: |e|=2\}$. We show that the polynomial $p_2$ is convex over the standard simplex or, equivalently, that its Hessian matrix is positive semidefinite over $\Delta_m$. Here,
 the Hessian matrix of $p_2$ is  the matrix indexed by $E$, renamed  $M$, with entries
\begin{equation}\label{eqM0}
M_{e,f}={1\over |e\cup f|}\quad \text{ for } e,f\in E.
\end{equation}
Consider  the matrices $A_2,A_3,A_4$ indexed by $E$, with entries
$$(A_s)_{e,f}=1 \ \text{ if } |e\cup f|=s, \quad (A_s)_{e,f}=0\ \text{ otherwise},\quad \text{ for } s=2,3,4.$$
Then, we have $A_2=I$ and $A_2+A_3+A_4=J$. Clearly we can express the matrix $M$ as a linear combination of these matrices:
\begin{equation}\label{eqM}
M={1\over 2}I+{1\over 3}A_3+{1\over 4}A_4= {1\over 4} I+{1\over 12}A_3+ {1\over 4} J
={1\over 12}I +{1\over 4}J +{1\over 12}(A_3+2I).
\end{equation}
We can now conclude that $M\succeq 0$ (and thus the polynomial $p_2$ is convex) in view of the next lemma, which claims that $A_3+2I\succeq 0$.

\begin{lemma}\label{lemJohnson}
Consider the  ${n\choose 2}\times n$ matrix $\Gamma_n$, with entries
$(\Gamma_n)_{e,i}= |e \cap \{i\}|$ for $e\in E$ and $i\in [n].$
Then  $A_3+2I = \Gamma_n\Gamma_n^T \succeq 0$.
\end{lemma}

\begin{proof}
Direct verification. \end{proof}

Note that the matrices $A_2=I, A_3,A_4$  generate the Bose-Mesner algebra of the Johnson scheme $J^n_2$, with length $n$  and weight $2$, and thus the matrix $M$ belongs to this Bose-Mesner algebra  (see \cite{DelsarteLevenshtein} for details on the Johnson scheme).
For  arbitrary degree $d\ge 3$ and edge size  $L= 2$ one could proceed to show that the Hessian matrix of $p_d$ is convex by using a similar symmetry reduction based on the Bose-Mesner algebra of the Johnson scheme $J^{p}_2$ for suitable values of $p$. However, for general edge size $L\ge 3$ we will need to use a richer algebra, namely the Terwilliger algebra of the Hamming cube. Hence we will treat in the rest of the section  the general case $d\ge 2$ and $L\ge 2$.

\subsection{Computing the Hessian matrix of $p_d$}\label{secHessian}

In this section  we indicate how to compute the Hessian matrix of the polynomial
\begin{equation}\label{eqpd}
p_d(x)=\sum_{(e_1,\ldots,e_d) \in E^d} c_{(e_1,\ldots,e_d)} x_{e_1}\cdots x_{e_d},
\end{equation}
where we set
\begin{equation}\label{eqc}
c_{(e_1,\ldots,e_d)}={1\over |e_1\cup \ldots \cup e_d|} \quad \text{ for } e_1,\ldots,e_d\in E
\end{equation}
and as before $E=\{e\subseteq V=[n]: |e|=L\}$ with $L\ge 2$.
We begin with getting the explicit coefficients of the polynomial $p_d$ expressed in the standard monomial basis. The basic fact we will now use is that the parameter $c_{(e_1,\ldots,e_d)}$ depends only on the set of distinct indices $e_i$ that are present in the tuple $(e_1,\ldots,e_d)\in E^d$ and not on their multiplicities.

To formalize this, recall $m=|E|$ and label the edges as $e_1,\ldots,e_m$ so that  $E=\{e_1,\ldots,e_m\}$. For a $d$-tuple $\underline e:=(e_{i_1},\ldots,e_{i_d})\in E^d$ with $i_1,\ldots,i_d\in [m]$, define the sequence $\alpha (\underline e)\in \oN^m$, where, for $\ell\in [m]$,  $\alpha(\underline e)_\ell $ is the number of indices among ${i_1}, \ldots, {i_d}$ that are equal to $\ell$. Then we have:
$$x_{e_{i_1}}\cdots x_{e_{i_d}} = x_{e_1}^{\alpha(\underline e)_1} \cdots x_{e_m}^{\alpha(\underline e)_m}=x^{\alpha(\underline e)}$$
and $|\alpha(\underline e)|=d$ so that $\alpha(\underline e)\in \oN^m_d$. This justifies the following definition.
For $\alpha\in \oN^m_d$, consider a $d$-tuple  $\underline e\in E^d$ such that $\alpha(\underline e)=\alpha$ and define
\begin{equation}\label{eqchat}
\widehat c_\alpha:= c_{\underline e}.
\end{equation}
 In this way we get $\widehat c=(\widehat c_\alpha)_{\alpha \in \oN^m_d}$ corresponding to the vector $c=(c_{\underline e})_{\underline e\in E^d}$
 in (\ref{eqc}). As an example, for $d=n=m=3$, if $\alpha =(1,0,2)$ then
$\widehat c_\alpha= c_{(e_1,e_3,e_3)}={1\over |e_1\cup e_3|}$. And, also for $\alpha=(2,0,1)$, we have
$\widehat c_\alpha= c_{(e_1,e_1,e_3)}={1\over |e_1\cup e_3|}$.

 We can now reformulate the polynomial $p_d$ in the (usual) monomial basis.

\begin{lemma}\label{lempd}
The polynomial $p_d$ from (\ref{eqpd}) can be reformulated as follows:
\begin{equation}\label{eqpd2}
p_d(x)=\sum_{\alpha\in \oN^m_d} \widehat c_\alpha {d!\over \alpha!}x^\alpha,
\end{equation}
setting $\alpha!=\alpha_1!\cdots \alpha_m!$ and where $\widehat c_\alpha$ is as defined in (\ref{eqchat}).
\end{lemma}

\begin{proof}
Using the definition of the coefficients $\widehat c_\alpha$, we  can rewrite $p_d$ as
$$p_d(x)=\sum_{\alpha\in \oN^m_d} \Big(\sum_{\underline e\in E^d: \alpha(\underline e)=\alpha} c_{\underline e}\Big) x^\alpha=
\sum_{\alpha\in \oN^m_d} \Big(\sum_{\underline e\in E^d: \alpha(\underline e)=\alpha} \widehat c_\alpha\Big) x^\alpha,
$$
which is equal to $ \sum_{\alpha\in \oN^m_d} \widehat c_\alpha {d!\over \alpha!} x^\alpha.$ Here, for this last equality, we use the monomial theorem, which claims
  \[
  \left(\sum_{i=1}^{m}x_i\right)^d = \sum_{\alpha\in\oN^m_d }\frac{d!}{\alpha!}x^\alpha,
  \]
  or, equivalently, that the number of $d$-tuples $\underline e\in E^d$ for which $\alpha(\underline e)=\alpha$ is equal to $d!/\alpha!$.
\end{proof}

We now proceed to compute the Hessian matrix of $p_d$.

\begin{lemma}\label{lemHpd}
The Hessian of the polynomial $p_d$ is the matrix
$$H(p_d)(x)= \Big({\partial^2 p_d(x)\over \partial x_{e_i} \partial x_{e_j}}\Big)_{i,j=1}^m =\sum_{\gamma\in \oN^m_{d-2}} {d!\over \gamma!}x^\gamma M_\gamma,
$$
where we set
\begin{equation}\label{eqMgamma}
M_\gamma =(\widehat c_{\gamma +u_i+u_j})_{i,j=1}^m
\end{equation}
and the vectors $u_1,\ldots,u_m\in \oR^m$ form the standard basis of $\oR^m$.
\end{lemma}

\begin{proof}
  The partial derivatives of $p_d$ are
  \[
  \frac{\partial p_d(x)}{\partial x_{e_i}}=\sum_{\alpha\in \oN^m_d: \alpha_i\ge 1 }\frac{d!}{(\alpha -u_i)!}\widehat c_{\alpha}x^{\alpha-u_i}
=  \sum_{\beta\in\oN^m_{d-1}}  \frac{d!}{\beta!}\widehat c_{\beta+u_i}x^{\beta}.
  \]
Similarly we see that
  \begin{align*}
    \frac{\partial^2 p(x)}{\partial x_{e_j}\partial x_{e_i}} & = \sum_{\beta\in\oN^m_{d-1}: \beta_j\ge 1}  \frac{d!}{(\beta -u_j)!}\widehat c_{\beta+u_i} x^{\beta-u_j}
    = \sum_{\gamma\in\oN^m_{d-2}} \widehat c_{\gamma+u_i+u_j} {d!\over \gamma!} x^\gamma.
  \end{align*}
  This concludes the proof.
\end{proof}

Hence, if we can show that the matrices $M_\gamma$ in (\ref{eqMgamma}) are all positive semidefinite then it  follows directly that the Hessian matrix of $p_d$ is positive semidefinite on the standard simplex.

We next observe how to further simplify the matrices $M_\gamma$. For $\gamma\in \oN^m$, define its support  as the set
$S_\gamma=\{e\in E: \gamma_e\ge 1\}$ and  let
$$W_\gamma=\bigcup_{e\in S_\gamma} e$$
denote the subset of elements of $V=[n]$ that are covered by some edge in $S_\gamma$. Then, for any $i,j\in [m]$, the support of $\gamma +u_i+u_j$ is the set $S_\gamma \cup\{e_i,e_j\}$ and we have
$$\widehat c_{\gamma+u_i+u_j}={1\over |W_\gamma \cup e_i \cup e_j|}.$$
Hence the matrix $M_\gamma$  depends only on the set $W_\gamma$ (and not on the specific choice of $\gamma$). This justifies defining the matrices
\begin{equation}\label{eqMW}
M_W= \Big({1\over |W\cup e\cup f|}\Big)_{e,f\in E}
\end{equation}
for any set $W\subseteq V=[n]$. 
Summarizing, we have shown:

\begin{lemma}\label{proppdconvex}
Assume that  the matrices $M_W$ from (\ref{eqMW}) are positive semidefinite for all $W\subseteq V$ with $|W|\ge L$ (if $d\ge 3$) and $|W|\le L(d-2).$
Then the polynomial $p_d$ is convex over the standard simplex.
\end{lemma}


If $d=2$ then there is only one matrix  to check, namely the matrix $M_\emptyset$ (for $W=\emptyset$). Note that the matrix $M_\emptyset$ coincides with the matrix in (\ref{eqM0}), so we already know it is positive semidefinite when $L=2$.
 However, if $d\ge 3$, then one needs to check all the matrices of the form $M_W$ in (\ref{eqMW}).
Now comes the last reduction, useful to link these matrices $M_W$ to the Terwiliger algebra,
which consists in removing duplicate rows and columns.
Set $p :=|W|$ and
  $U:=V\setminus W$, so that  $|U|=n-p$. 
In addition let
\begin{equation}\label{eqF}
\FmL:=\{e\subseteq U: |e|\le L\}
\end{equation}
denote the collection of subsets of $U$ with size at most $L$.
Now we consider the matrix $M_p$, which is  indexed by $\FmL$, with entries
\begin{equation}\label{eqMp}
(M_p)_{e,f}={1\over p+|e\cup f|}\quad \text{ for } e,f\in \FmL.
\end{equation}
Note that for $p=0$ the matrix $M_0$ coincides with the matrix $M_\emptyset$ in (\ref{eqMW}) (and  with the matrix in  (\ref{eqM0})).
The next lemma links the matrices $M_W$ and $M_p$.

\begin{lemma}\label{symmRedu}
	Let $L \geq 2$ and $d \geq 2$. Consider the matrces $M_W$ in (\ref{eqMW}) and $M_p$ in (\ref{eqMp}).
	The following assertions are equivalent:
	\begin{itemize}
		\item[(i)]  $M_W\succeq 0$ for all $W=e_1\cup \ldots \cup e_{d-2}$ with $e_1,\ldots e_{d-2}\in \EnL$.
		\item[(ii)]  $M_p\succeq 0$ for all $p\le L(d-2)$ such that  $p\ge L$ if $d\ge 3$.
	\end{itemize}
\end{lemma}

\begin{proof}
Let $W=e_1\cup \ldots \cup e_{d-2}$, where $e_1,\ldots,e_{d-2}\in E$.
Consider the partition of the set $E$ into
$E=\cup_{i=0}^L E_i$, where $E_i=\{e\in E: |e\setminus W|=i\}.$
With respect to this partition of its index set, the matrix $M_W$
has the following block-form:
	\[
	M_W =
	\left(
	\begin{array}{c|c|c|c}
	M^{0,0}_{W} & M^{0,1}_{W} & \cdots & M^{0,L}_{W}  \\
	\hline
	M^{1,0}_{W} & M^{1,1}_{W} & \cdots & M^{1,L}_{W}  \\
	\hline
	\vdots  & \vdots     & \ddots & \vdots  \\
	\hline
	M^{L,0}_{W} & M^{L,1}_{W} & \cdots & M^{L,L}_{W} \\
	\end{array}
	\right),
	\]
	where the block $M^{i,j}_{W}$ has its rows indexed by $E_i$ and its columns by $E_j$.
Note that, if two edges $e,e'\in E$ satisfy $e\setminus W = e'\setminus W$, then the two rows of $M_W$ indexed  by $e$ and $e'$ coincide:
for any $f\in E$ we have
		$$(M^{i,j}_{W})_{e,f} = \frac{1}{|W|+|(e\cup f) \setminus W|} = \frac{1}{|W|+|(e'\cup f) \setminus W|} =  (M^{i,j}_{W})_{e',f}.$$
	In fact, after removing these duplicate rows (and columns) and keeping only one copy for each subset of $U=V\setminus W$,  we obtain the matrix
	
	\[
	\left(
	\begin{array}{c|c|c|c}
	M^{0,0}_{p} & M^{0,1}_{p} & \cdots & M^{0,L}_{p}  \\
	\hline
	M^{1,0}_{p} & M^{1,1}_{p} & \cdots & M^{1,L}_{p}  \\
	\hline
	\vdots  & \vdots     & \ddots & \vdots  \\
	\hline
	M^{L,0}_{p} & M^{L,1}_{p} & \cdots & M^{L,L}_{p} \\
	\end{array}
	\right),
	\]
	which coincides with the matrix $M_p$ in (\ref{eqMp}). Indeed, the above matrix is indexed by the set $F$ in (\ref{eqF}) and its block-form is with respect to the partition $F=\cup_{i=0}^L F_i$, where $F_i=\{e\subseteq U: |e|=i\}$.
	So the block $M^{i,j}_{p}$ has its rows indexed by $F_i$, its columns indexed by $F_j$, and its entries are
	\begin{equation}\label{eqMp1}
	(M^{i,j}_{p})_{e,f} ={1\over p+|e\cup f|}= {1\over p+i+j-|e\cap f|} \ \text{ for } e\in F_i, f\in F_j.
	\end{equation}
	As the matrices  $M_p$ arise from $M_W$ by removing its duplicates rows and columns  it is clear that the matrices $M_W$ are positive semidefinite if and only if the same holds for the matrices $M_p$. This concludes the proof.
\end{proof}

\subsection{The general case $d\ge 2$ and $L\ge 2$}\label{secdge3Lge3}

In Section \ref{secPrelimTer} we gave preliminary results on the Terwilliger algebra, which we will now use to prove that the matrices $M_p$ in (\ref{eqMp}) are positive semidefinite. Fix $0\le p\le n$ and consider the matrix $M_p$ in (\ref{eqMp}) (with blocks as in (\ref{eqMp1})). We start with observing  that $M_p$ belongs to  the Terwilliger algebra  $\mathcal{A}_{n-p}$. This is clear since relation (\ref{eqMp1}) provides the explicit  correspondence between the blocks $M^{i,j}_{p}$ of $M_p$ and the generating matrices $D_{i,j}^t$ of the algebra $\mathcal A_{n-p}$:
$$
M_p =  \sum_{i = 0}^{L}\sum_{j = 0}^{L}\sum_{t = 0}^{\uma} \frac{1}{p +i + j - t } D_{i,j}^t
=  \sum_{i = 0}^{L}\sum_{j = 0}^{L}\sum_{t = 0}^{\uma}x_{i,j}^t D_{i,j}^t,
$$
after setting
\begin{equation}\label{eqx}
x_{i,j}^t=  \frac{1}{p +i + j - t }.
\end{equation}
Let $B_k$ be the corresponding matrices from (\ref{eqBk}). Then, in view of  Theorem~\ref{schrter}, we know that $M_p\succeq 0$ if and only if $B_k\succeq 0$ for all $0\le k\le \lceil n/2\rceil.$



In what follows $p,k$ are fixed integers. We now proceed to show that $B_k\succeq 0$.

\medskip
To simplify the notation we introduce the following parameters
\begin{align*}
  \comba i) &:= {n-p-2k\choose i-k}^{-\frac{1}{2}}, &
  \combb i)& := {n-p-k-u \choose i-u},&
  \combc & := {n-p-2k \choose n-p-k-u}.
\end{align*}
for any  integers $i,u$. Note that we may omit the obvious bounding conditions on $i$ and $u$ since the corresponding parameters are zero if these conditions  are not satisfied; for instance, $a(i)=0$ if $i<k$ and $b(u,i)=0$ if $u>i$.
Then we have
\begin{equation}\label{eqBk}
	B_k = \left( \comba i) \comba j) \sum_{t=0}^{min\{i,j\}}    \beta_{i,j,k}^{t}x_{i,j}^{t}   \right)_{i,j = k}^{n-p-k}
\end{equation}
and
\begin{equation}\label{eqbeta}
    \beta_{i,j,k}^{t} := \sum_{u = 0}^{n-p} (-1)^{u-t} {u \choose t} \combc \combb i) \combb j).
\end{equation}
We now give an integral reformulation for the entries of the matrix $B_k$ from (\ref{eqBk}).  It is  based on the fact that
\begin{equation}\label{eqint}
{1\over i} =\int_0^1 z^{i-1} dz\quad \text{ for any integer } i\ge 1,
\end{equation}
which permits to give an integral reformulation for the scalars $x_{i,j}^t$ in (\ref{eqx}). This simple but powerful fact  will be very useful  to show $B_k\succeq 0$.

\begin{lemma}\label{lemint1}
We have $$
\sum_{t=0}^{min\{i,j\}}\beta_{i,j,k}^{t}x_{i,j}^{t} = \sum_{u = 0}^{\uma} \combc \combb i)\combb j) \int_0^1  \g z^{i+j}  dz,
$$		
where we define the function  $\g = z^{p-1} (\frac{1-z}{z})^u$ for $z\in (0,1]$. 
\end{lemma}

\begin{proof}
First we use the expressions of $\beta_{i,j,k}^t$ in (\ref{eqbeta}) and of $x_{i,j}^t$ in (\ref{eqx}) and we exchange the summations in $t$ and $u$ to obtain
\begin{equation}\label{eq1}
\sum_{t=0}^{\uma}\beta_{i,j,k}^{t}x_{i,j}^{t}
=  \sum_{u = 0}^{\uma} \Bigg{(} \sum_{t=0}^{u} \frac{1}{p +i + j - t } (-1)^{u-t} {u \choose t} \Bigg{)} \combc \combb i)\combb j).
\end{equation}
Now we use (\ref{eqint}), which gives the following integral representation
$$
 \frac{1}{p +i + j - t } = \int_0^1z^{p+i+j-t-1}dz.
$$
Using this integral representation (and the binomial theorem for the equality marked (*) below) we can reformulate  the inner summation appearing in (\ref{eq1}) as follows:
\begin{align*}
  \sum_{t=0}^{u} \frac{1}{p +i + j - t } (-1)^{u-t} {u \choose t} & =\sum_{t=0}^{u}(-1)^{u-t} {u \choose t}      \int_0^1z^{p+i+j-t-1}dz \\
   & = \int_0^1  z^{p+i+j-1} (-1)^{u}\left( \sum_{t=0}^{u} \left(-\frac{1}{z}\right)^{t} {u \choose t} \right) dz \\
   & \stackrel{(*)}{=} \int_0^1z^{p+i+j-1}(-1)^{u} \left(1-\frac{1}{z}\right)^u dz\\
   & = \int_0^1z^{p+i+j-1}(-1)^{u} \left(\frac{z-1}{z}\right)^u dz\\
   &= \int_0^1  z^{p-1} \left(\frac{1-z}{z}\right)^u z^{i+j} dz.
\end{align*}
This concludes the proof.
\end{proof}

We can now  proceed to show  that the matrices $B_k$ in (\ref{eqBk}) are positive semidefinite.

\begin{lemma}
We have 	$
	B_k \succeq 0.
	$		
\end{lemma}

\begin{proof}
We use  Lemma \ref{lemint1} to reformulate the  matrix $B_k$. First, note that in the result of Lemma \ref{lemint1}, since $b(u,i)b(u,j)=0$ if $u>\min\{i,j\}$, we may replace the summation on $u$ from $0\le u\le \min\{i,j\}$ to $0\le u\le n-p$. This implies:
\begin{align*}
  B_k& = \Big( \comba i)\comba j)\sum_{t= 0}^{n-p}  \beta_{i,j,k}^{t}x_{i,j}^{t}   \Big)_{i,j = k}^{n-p-k} \\
   & = \int_0^1 \Big( \sum_{u=0} ^{n-p} \g \combc \underbrace{(z^i \comba i)\combb i))}_{=:h(u,z,i)}
   \underbrace{(z^j \comba j)\combb j)) }_{=:h(u,z,j)} \Big)_{i,j = k}^{n-p-k} dz \\
   &= \sum_{u= 0}^{n-p} \int_0^1  \g \combc \underbrace{ \Big(h(u,z,i)h(u,z,j) \Big)_{i,j = k}^{n-p-k}}_{=: H(u,z,k)} dz \\
   &= \sum_{u\geq 0} \int_0^1  \underbrace{\g \combc}_{\geq 0} \underbrace{H(u,z,k)}_{\succeq 0} dz\succeq 0.
\end{align*}
 Here we used the fact that, for any $u\in [0,n-p]$, the function $g(u,z)$ is nonnegative on $(0,1]$ and that  the matrix $H(u,z,k)$ is positive semidefinite for any $z\in [0,1]$ since it is the outerproduct of the vector $h(u,z,i)$ with itself.
\end{proof}

Therefore we have shown that the matrices $B_k$ are positive semidefinite and thus the following result.

\begin{corollary}\label{corMp}
The matrices $M_p$ from (\ref{eqMp}) are positive semidefinite for all $0\le p\le n$.
\end{corollary}

In view of Lemmas \ref{proppdconvex} and \ref{symmRedu} we can conclude  that the polynomial  $p_d$ is convex on $\Delta_m$,
which concludes the proof of Theorem~\ref{theomainpd}.

\section{Investigating the polynomials $f_d$}\label{secprooffd}

Here we consider the second class  of  polynomials $f_d$ from (\ref{eqfd0}), namely
$$ f_d(x) =\sum_{(e_1,\ldots,e_d)\in E^d}\prod_{i=1}^{d}\frac{x_{e_i}}{\vert e_1\cup\ldots\cup e_i\vert}.$$
We address Question \ref{question}, which asks  whether $f_d$ attains its minimum value on the simplex $\Delta_m$ at the barycenter of  $\Delta_m$.
Here too this question has  a positive answer if one can show that $f_d$ is convex over $\Delta_m$. This follows since the analogue of Lemma \ref{lemcondition} extends easily  for the polynomial $f_d$. We conjecture that convexity holds in general.

\begin{conjecture}\label{conjfd}
For any integers $n,L,d\ge 2$ the polynomial $f_d$ is convex over the simplex $\Delta_m$.
\end{conjecture}

For degree $d=2$, we have $f_2={1\over L} p_2$ and thus we know from Theorem~\ref{theomainpd} that $f_2$ is convex.
We  will  prove in Section   \ref{secfd3}  that Conjecture \ref{conjfd} holds   for  degree $d=3$ and edge size $L=2$ and, in Section \ref{secnumeric} and Appendix A,  we will give a range of values for $(n,L,d)$ that were numerically tested and support Conjecture \ref{conjfd}.

In what follows we begin  in Section \ref{secHfd} with giving a polynomial matrix decomposition for the Hessian of $f_d$ and then a  recursive reformulation for it, making also apparent some links to the Hessian of $p_d$. From this we see that convexity of $f_d$ follows if we can show that a family of well-structured matrices $Q_\gamma$ are positive semidefinite (see Lemmas \ref{lemHfd} and \ref{lembahat}). We can complete this task in the case $d=3$ and $L=2$ (see Section \ref{secfd3}).
However, understanding the general case is technically involved and would require developing new tools for exploiting the symmetry structure present in the matrices $Q_\gamma$ (which is now not captured by the Terwilliger algebra). This goes beyond the scope of this paper and we leave it for further research.

\subsection{Computing the Hessian of $f_d$}\label{secHfd}


We begin with expressing the polynomial $f_d$ in the standard monomial basis:
\begin{equation}\label{eqfd1}
 f_d(x)  =\sum_{\alpha\in \mathbb{N}^m_d}x^\alpha
\sum_{\substack{\underline{e}=(e_1,\ldots,e_d)\in E^d\\ \alpha(\underline{e})=\alpha}}\prod_{i=1}^{d}\frac{1}{\vert e_1\cup\ldots\cup e_i\vert}
=\sum_{\alpha\in \mathbb{N}^m_d}b_\alpha x^\alpha,
\end{equation}
where we set
\begin{equation}\label{eqbalpha}
b_\alpha= \sum_{\substack{\underline{e}=(e_1,\ldots,e_d)\in E^d\\ \alpha(\underline{e})=\alpha}}\prod_{i=1}^{d}\frac{1}{\vert e_1\cup\ldots\cup e_i\vert}.
\end{equation}
Next we  compute the Hessian of $f_d$ and we give a matrix polynomial reformulation for it.

\begin{lemma}\label{lemHfd}
The Hessian  of the polynomial $f_d$ is given by
$${\partial^2 f(x)\over \partial x_{e_i}\partial x_{e_j}}=
\left\{
\begin{array}{ll}
 \sum_{\gamma\in \oN^m_{d-2}} (\gamma_i+1)(\gamma_j+1) x^\gamma b_{\gamma +u_i+u_j} \quad \text{ if } i\ne j\\
 \sum_{\gamma\in \oN^m_{d-2}} (\gamma_i+1)(\gamma_i+2) x^\gamma b_{\gamma +2u_i} \quad \text{ if } i=j
 \end{array}
 \right.
 $$
 where as before $u_1,\ldots,u_m$ form  the standard basis of $\oR^m$.  In other words,
 $$H(f_d) (x)=\sum_{\gamma \in \oN^m_{d-2}} x^\gamma Q_\gamma,$$
 where we define the symmetric $m\times m$ matrix $Q_\gamma$ with entries
\begin{equation}\label{eqQgamma}
(Q_\gamma)_{ij}= (\gamma_i+1)(\gamma_j+1) b_{\gamma+u_i+u_j} \text{ if } i\ne j,\ \
(Q_{\gamma})_{ii}= (\gamma_i+1)(\gamma_i+2) b_{\gamma+2u_i}
\end{equation}
for $i,j\in [m]$ and $\gamma\in \oN^m_{d-2}$.
Hence, $H(f_d)\succeq 0$ if  $Q_\gamma \succeq 0$ for all $\gamma\in\oN^m_{d-2}$.
  \end{lemma}

\begin{proof}
Direct verification.
\end{proof}

We now give a recursive reformulation for the coefficients of the polynomial $f_d$ and for its Hessian matrix, that may  possibly be helpful for a proof by induction.
Recall the definition of the coefficients $b_\alpha$ of $f_d$ in (\ref{eqbalpha}). Fix $\alpha\in \oN^m_d$.
 There are $d!\over \alpha !$ distinct tuples
 $\underline e$ such that $\alpha(\underline e)=\alpha$. For any  such sequence
$\underline e=(e_{i_1},\ldots, e_{i_d})$ with $i_1,\ldots,i_d\in [m]$, $\alpha=\alpha(\underline e)$    means that, for any $\ell\in [m]$,   $\alpha_\ell $ is the number of occurrences of $\ell$ within the multiset $\{i_1,\ldots,i_d\}$; so $\alpha_\ell \ge 1$ if $\ell \in \{i_1,\ldots,i_d\}$ and $\alpha_\ell=0$ if $\ell\not\in \{i_1,\ldots,i_d\}$.
For instance, for $\underline e=(e_1,e_2,e_3,e_2,e_1)$, $d=5$, $m=4$, we have
$(i_1,\ldots,i_5)= (1,2,3,2,1)$ and
$\alpha(\underline e) = (2,2,1,0)$.


To reformulate $b_\alpha$ we exploit the fact that $b_\alpha$ enjoys some  invariance property  under permutations of $[d]$, namely
\begin{align}
  b_\alpha & = \sum_{\substack{\underline{e}=(e_{i_1},\ldots,e_{i_d})\in E^d:\\ \alpha(\underline{e})=\alpha}}\prod_{k=1}^{d}\frac{1}{\vert e_{i_1}\cup\ldots\cup e_{i_k}\vert} \\
   & = \frac{1}{d!}\sum_{\sigma\in \Sym(d)}\sum_{\substack{\underline{e}=(e_{i_1},\ldots,e_{i_d})\in E^d:\\ \alpha(\underline{e})=\alpha}}\prod_{k=1}^{d}\frac{1}{\vert e_{i_{\sigma(1)}}\cup\ldots\cup e_{i_{\sigma(k)}}\vert}\\
   & = \frac{1}{d!}\sum_{\substack{\underline{e}=(e_{i_1},\ldots,e_{i_d})\in E^d\\ \alpha(\underline{e})=\alpha}}\underbrace{\sum_{\sigma\in \Sym(d)}\prod_{k=1}^{d}\frac{1}{\vert e_{i_{\sigma(1)}}\cup\ldots\cup e_{i_{\sigma(k)}}\vert}}_{=:S}.\label{eqb3}
   \end{align}
Observe that the inner summation $S$ in (\ref{eqb3}) does not depend on the choice of the sequence $\underline e$ such that $\alpha(\underline e)=\alpha$; thus we may consider it fixed, denoted as $(e_{i_1},\ldots,e_{i_d})$. Since there are $d!\over \alpha!$ possible choices for selecting this sequence,  using  relation (\ref{eqb3}) we can reformulate $b_\alpha$ as follows:
\begin{align*}
b_\alpha   &=
{1\over d!}{d!\over \alpha!}   \sum_{\sigma\in \Sym(d)}\prod_{k=1}^{d}\frac{1}{\vert e_{i_{\sigma(1)}}\cup\ldots\cup e_{i_{\sigma(k)}}\vert}
    = {1\over \alpha!}  \sum_{\sigma\in \Sym(d)}\prod_{k=1}^{d}\frac{1}{\vert e_{i_{\sigma(1)}}\cup\ldots\cup e_{i_{\sigma(k)}}\vert}.
      \end{align*}
   Next we pull out the factor ${1\over |e_{i_1}\cup\ldots \cup e_{i_d}|}=\widehat c_\alpha$ which occurs for $k=d$ and get
   \begin{align*}
b_\alpha
   &= \frac{\widehat c_\alpha}{\alpha!} \sum_{r=1}^{d}\sum_{\sigma\in \Sym(d):\sigma(d)=r}\prod_{k=1}^{d-1}\frac{1}{\vert e_{i_{\sigma(1)}}\cup\ldots\cup e_{i_{\sigma(k)}}\vert}\\
   &=\frac{\widehat c_\alpha}{\alpha!} \sum_{r=1}^{d}    b_{\alpha-u_{i_r}} (\alpha-u_{i_r})!\\
   & = \widehat c_\alpha \sum_{r=1}^d {b_{\alpha-u_{i_r}}\over \alpha_{i_r}}\\
   & \stackrel{(*)}{=} \widehat c_\alpha \sum_{k\in [m]: \alpha_k\ge 1} b_{\alpha -u_k}. 
\end{align*}
Here, in the last equality marked (*),  we use the fact that 
$\alpha_k$ of the elements in  the multiset $\{i_1,\ldots,i_d\}$  are equal to $k$.
Summarizing we have shown:

\begin{lemma}\label{lemba0}
For any   $\alpha\in \oN^m_d$ 
 we have
$$b_\alpha=
 \widehat c_\alpha \sum _{k\in [m]:\alpha_k\ge 1} b_{\alpha-u_k}.$$
\end{lemma}

We now proceed to give a recursive reformulation for the matrices $Q_\gamma$ in (\ref{eqQgamma}). First we reformulate them using the scaled parameters
\begin{equation}\label{eqbalphahat}
\widehat b_\alpha:=  \alpha! \ b_\alpha,
\end{equation}
which satisfy the recursive relation:
\begin{equation}\label{eqbalphahatrec}
\widehat b_\alpha=\widehat c_\alpha \sum_{k: \alpha_k\ge 1} \alpha_k  \widehat b_{\alpha-u_k}.
\end{equation}
Indeed, by Lemma \ref{lemba0} we have
$$\widehat b_\alpha= \alpha!\ b_\alpha=\alpha!\  \widehat c_\alpha \sum_{k:\alpha_k\ge 1} b_{\alpha-u_k}=
 \alpha!\ \widehat c_\alpha  \sum_{k:\alpha_k\ge 1} {\widehat b_{\alpha-u_k} \over \alpha-u_k!}
= \widehat c_\alpha \sum_{k: \alpha_k\ge 1}\alpha_k  \widehat b_{\alpha-u_k}.$$

\begin{lemma} \label{eqQgammahat}
For any $\gamma\in\oN^m_{d-2}$ we have $Q_\gamma={1\over \gamma!} (\widehat b_{\gamma +u_i+u_j})_{i,j=1}^m.$
\end{lemma}

\begin{proof}
Direct verification: for $i\ne j$ we have $(Q_\gamma)_{ij}=(\gamma_i+1)(\gamma_j+1)b_{\gamma+u_i+u_j}= \widehat b_{\gamma+u_i+u_j}(\gamma_i+1)(\gamma_j+1)/(\gamma+u_i+u_j)!= b_{\gamma+u_i+u_j}/\gamma!$ and, for $i=j$, we have
$(Q_\gamma)_{ii}=(\gamma_i+1)(\gamma_i+2)b_{\gamma+2u_i} =\widehat b_{\gamma+2u_i} (\gamma_i+1)(\gamma_i+2)/(\gamma+2u_i)!= \widehat b_{\gamma+2u_i} /\gamma!$.
\end{proof}

\begin{lemma}\label{lembahat}
For $d\ge 3$ and   $\gamma\in\oN^m_{d-2}$ we have
$$Q_\gamma= \underbrace{(\widehat c_{\gamma +u_i+u_j})_{i,j=1}^m}_{M_\gamma}
 \circ
 \Big(
 \sum_{k\in [m]: \gamma_k\ge 1} Q_{\gamma-u_k}
 +
\underbrace{ {1\over \gamma!}(\widehat b_{\gamma+u_i} +\widehat b_{\gamma+u_j})_{i,j=1}^m}_{=:R_\gamma}
 \Big) $$
 $$= M_{\gamma}\circ \Big( \sum_{k\in [m]: \gamma_k\ge 1} Q_{\gamma-u_k}  +R_{\gamma}\Big),$$
 where the matrices $M_\gamma$ were introduced in (\ref{eqMgamma}).
 \end{lemma}

\begin{proof}
Combining Lemmas \ref{lemba0} and \ref{lembahat}  we obtain
\begin{align*}
(Q_\gamma)_{ij}& ={1\over \gamma!} \widehat b_{\gamma+u_i+u_j}= {1\over\gamma!}\widehat c_{\gamma+u_i+u_j} \sum_{k: (\gamma+u_i+u_j)_k\ge 1} \widehat b_{\gamma+u_i+u_j-u_k} (\gamma+u_i+u_j)_k\\
& ={1\over \gamma!}\widehat c_{\gamma+u_i+u_j} \Big(\sum_{k\ne i,j: \gamma_k\ge 1} \widehat b_{\gamma+u_i+u_j-u_k} \gamma_k
+\widehat b_{\gamma+u_j} (\gamma_i+1)+\widehat b_{\gamma+u_j}(\gamma_i+1)\Big)\\
&={1\over \gamma!} \widehat c_{\gamma+u_i+u_j} \Big(\sum_{k: \gamma_k\ge 1} \widehat b_{\gamma-u_k+u_i+u_j} \gamma_k
+ \widehat b_{\gamma+u_i}+\widehat b_{\gamma+u_j}\Big)\\
&= \widehat c_{\gamma+u_i+u_j} \Big(\sum_{k: \gamma_k\ge 1}{ \widehat b_{\gamma-u_k+u_i+u_j} \over (\gamma-u_k)!}
+ {1\over \gamma!} (\widehat b_{\gamma+u_i}+\widehat b_{\gamma+u_j})\Big)\\
&=  \widehat c_{\gamma+u_i+u_j} \Big(\sum_{k: \gamma_k\ge 1} (Q_{\gamma-u_k})_{ij} +{1\over \gamma!} (\widehat b_{\gamma+u_i}+\widehat b_{\gamma+u_j})\Big),
\end{align*}
which shows the claim.
\end{proof}

\subsection{The polynomial $f_d$ in the case d = 3, L = 2}\label{secfd3}

 Here  we show that the polynomial $f_d$ is convex in the case $d=3$ and $L=2$.  In view of Lemma \ref{lemHfd} it suffices to show that the matrix $Q_\gamma$ is positive semidefinite for any $\gamma \in \oN^m_1$. Up to symmetry it suffices to show that $M_\gamma\succeq 0$ for $\gamma=u_1$.
  In view of Lemma \ref{lembahat} we have
  $$Q_{u_1}= \underbrace{ (\widehat c_{u_1+u_i+u_j})_{i,j=1}^m}_{=M_{u_1}}
   \circ
  (Q_0+ \underbrace{(\widehat b_{u_1+u_i}+\widehat b_{u_1+u_j})_{i,j=1}^m }_{=R_{u_1}}).
  $$
    where 
 $M_{u_1}=M_2$ is as defined in relation (\ref{eqMp}).
We have seen in  the previous section  that the matrix $M_{u_1}$ is positive semidefinite (see Corollary \ref{corMp}). Hence it suffices now to show that $Q_0+R_{u_1}\succeq 0$.
By definition, the entries of $Q_0$ (case $\gamma=0$) are
$$(Q_0)_{ii}= 2b_{2u_i}= {2\over L},\ \ (Q_0)_{ij} = b_{u_i+u_j}= {2\over |e_i\cup e_j|} \ \text{ for } i\ne j \in [m].$$
Moreover, $\widehat b_{2u_1}= 2b_{2u_1}={2\over L}$ and $\widehat b_{u_1+u_i}= b_{u_1+u_i}={2\over |e_1\cup e_i|}$ for $i\ge 2$.
Using this we obtain that
\begin{equation*}\label{eqQ0R}
Q_0+R_{u_1}=2\cdot  \Big( {1\over |e_1\cup e_j|}+{1\over |e_i\cup e_j|}+{1\over |e_1\cup e_i|}\Big)_{i,j=1}^m =: 2B,
\end{equation*}
where we define the matrix $B$ as
\begin{equation}\label{eqBf3}
 B:= \Big( \frac{1}{|e \cup f|} + \frac{1}{|e_1 \cup e|} + \frac{1}{|e_1 \cup f|}        \Big)_{e,f \in E}.
\end{equation}
The main result of this section is the next lemma, which implies that the polynomial $f_3$ is convex for $L=2$ and thus settles Conjecture \ref{conjfd} for the case $d=3,L=2$.

\begin{lemma}\label{lemmainf3}
Assume $L=2$. The matrix $B$ in (\ref{eqBf3}) is positive semidefinite.
\end{lemma}

Before proceeding to the proof, let us make a few observations.
Note that $B$ can be decomposed as $B=M_0+R,$ where $M_0=M_\emptyset$ has been shown earlier to be positive semidefinite (recall Corollary \ref{corMp}, or note that $M_0$ is the matrix $M$ from (\ref{eqM0})
as we are in the case $L=2$). On the other hand, the matrix $R$ is not positive semidefinite. In fact, $R$ has rank 2 and it has a negative eigenvalue. One can infer from the results in Section \ref{secd2L2} that $\lambda_{\min}(M_0) =1/12$, while one can check that $\lambda_{\min}(R)<-1/12 \ =-0.0833...$ when $n \ge 6$ (see Table \ref{tabled3L2}).
Hence in general one cannot  argue that $B\succeq 0$ by simply looking at the smallest eigenvalues of its summands $M_0$ and $R$.

\medskip
In the rest of the section we prove Lemma \ref{lemmainf3}.
To fix ideas let $e$ be the edge $e_1 = \{1,2\}$ and to simplify notation set $p=n-2$ and $q={n-2\choose 2}$.
Then the index set of $B$ can be partitioned into $\{e_1\} \cup I_1\cup I_2 \cup I_0$, setting $I_k= \{\{k,i\}: 3\le i\le n\}$ for $k=1,2$, and
$I_0 = \{\{i,j\}: 3\le i<j\le n\}.$ So $|I_1|=|I_2|=p$ and $|I_0|=q$.
With respect to this partition one can verify that  the matrix $B$ has the following block-form:


\begingroup
\fontsize{13pt}{25pt}\selectfont
$$
B =
\begin{blockarray}{ccccc}
&e_1  & I_1 & I_2 & I_0\\
\begin{block}{c(c|c|c|c)}
e_1         & \frac{3}{2} 	  & \frac{7}{6} J_{1,p}   & \frac{7}{6} J_{1,p}                &J_{1,q} \\
 \cline{2-5}
I_1     & \frac{7}{6} J_{p,1}	  & J_{p} + \frac{1}{6} I_{p}    & \frac{11}{12}J_{p}  + \frac{1}{12} I_{p}  &  \frac{5}{6}J_{p,q} +\frac{1}{12}\Gamma^T  \\
 \cline{2-5}
I_2    & \frac{7}{6} J_{p,1}   & \frac{11}{12}J_{p}  +\frac{1}{12} I_{p}  & J_{p}  + \frac{1}{6} I_{p}    & \frac{5}{6}J_{p,q} +\frac{1}{12}\Gamma^T   \\
 \cline{2-5}
I_0 & J_{q,1}      	  & \frac{5}{6}J_{q,p} +\frac{1}{12}\Gamma & \frac{5}{6}J_{q,p}+\frac{1}{12}\Gamma 			 & M + \frac{1}{2}J_{ q  } 		  \\
\end{block}
\end{blockarray} \ .
$$
\endgroup
Here 
 the matrix $M$ is the matrix from (\ref{eqM0}) (replacing $n$ by $p=n-2$), i.e., 
$$
M =  \frac{1}{12}I_{q} + \frac{1}{4}J_{q}+ \frac{1}{12}\Gamma \Gamma^T,
$$
where $\Gamma=\Gamma_{p}$ is the $ {p \choose 2} \times p$
matrix from Lemma~\ref{lemJohnson},
whose $(f,i)$th entry   is $ |\{i\} \cap f|$.

\ignore{\begin{proof}
One can compute the $(e,f)$th entry of $B$  by direct case checking:
\begin{itemize}
	\item[1.] For $e=f=e_1$:  $B_{e,f} = \frac{3}{|\{1,2\}|}   = \frac{3}{2}$
	\item[2.] For $e=e_1,f = \{1,j\}\text{ or }\{2,j\}$:   $B_{e,f} = \frac{2}{|\{1,2,j\}|}   + \frac{1}{|\{1,2\}|} = \frac{7}{6}$
	\item[3.] For $e=e_1,f = \{i,j\}$:   $B_{e,f}= \frac{1}{|\{1,2\}|} + \frac{2}{|\{1,2,j_1,j_2\}|} = 1$
	\item[4.] For $e=\{1,i\},f = \{1,j\}$:   $B_{e,f} =  \frac{1}{|\{1,i,j\}|}  + \frac{1}{|\{1,2,i\}|} + \frac{1}{|\{1,2,j\}|} =  \frac{1}{6}\delta_{ij} + 1 $
	\item[5.] For $e=\{1,i\},f = \{2,j\}$:  $B_{e,f}=  \frac{1}{|\{1,2,i,j\}|}  + \frac{1}{|\{1,2,i\}|} + \frac{1}{|\{1,2,j\}|} = 1 - \frac{1}{12}(1-\delta_{ij})$
	\item[6.] For $e=\{1,i\},f = \{j_1,j_2\}$:  $B_{e,f} = \frac{1}{|\{1,i,j_1,j_2\}|}  + \frac{1}{|\{1,2,i\}|} + \frac{1}{|\{1,2,j_1,j_2\}|} =\frac{5}{6} +\frac{1}{12}(|i \cap f|)$
	\item[7.] For $e=\{i_1,i_2\},f = \{j_1,j_2\}$:  $B_{e,f} = \frac{1}{|e \cup f |} + \frac{1}{2}$.
\end{itemize}
\end{proof}	
}


We now proceed to show that the matrix $B$ is positive semidefinite. Note that its lower right diagonal block indexed by the set $I_0$ is positive semidefinite (since $M\succeq 0$).
Our strategy is now to `eliminate' the borders indexed by the sets $\{e_1\}$, $I_1$ and $I_2$ successively, one by one, by taking Schur complements, until reaching a final matrix whose positive semidefiniteness can be seen  directly.
To do the Schur complement operations we will need to invert matrices of the form $aI+bJ$. The next lemma indicates how to do that,
its proof is straightforward and thus omitted.

\begin{lemma}\label{leminverse}
For $a,b \in \oR$ such that $a+pb\ne 0$, the matrix $aI_p+bJ_p$ is nonsingular with inverse
$$
(aI_p+ bJ_p)^{-1} = \frac{1}{a}\Big(I_p -\frac{b}{pb+a}\Big)J_p.
$$
\end{lemma}

\medskip\noindent
{\bf Step 1:} We take a first Schur complement with respect to the upper left corner of $B$ (indexed by $e_1$) and call $\widetilde B_1$ the resulting matrix, which reads

\begingroup
\fontsize{10pt}{18pt}\selectfont
$$
\begin{blockarray}{ccc}
\begin{block}{(c|c|c)}
 J_{p} + \frac{1}{6} I_{p}    & \frac{11}{12}J_{p}  + \frac{1}{12} I_{p}  &  \frac{5}{6}J_{p,q} +\frac{1}{12}\Gamma^T  \\
  \cline{1-3}
 \frac{11}{12}J_{p}  +\frac{1}{12} I_{p}  & J_{p}  + \frac{1}{6}I_p & \frac{5}{6}J_{p,q} +\frac{1}{12}\Gamma^T   \\
 \cline{1-3}
\frac{5}{6}J_{q,p} +\frac{1}{12}\Gamma & \frac{5}{6}J_{q,p} +\frac{1}{12}\Gamma 		 &  \frac{1}{12}I_q  + \frac{1}{12}\Gamma \Gamma^T + \frac{3}{4}J_{ q  } 		  \\
\end{block}
\end{blockarray}
$$
\vspace*{-0.4cm}
$$
-
\frac{2}{3}\begin{pmatrix}
\frac{7}{6} J_{p,1}   \\
 \frac{7}{6} J_{p,1} \\
J_{q,1}
\end{pmatrix}
\begin{pmatrix}
\frac{7}{6} J_{1,p}   & \frac{7}{6}J_{1,p}       & J_{1,q}\\
\end{pmatrix}
$$
\endgroup
\vspace*{-0.4cm}
 \begingroup
\fontsize{10pt}{18pt}\selectfont
 \begin{gather*}
=
\begin{blockarray}{ccc}
\begin{block}{(c|c|c)}
{ 5\over 54}  J_{p} + \frac{1}{6} I_{p}    & {1\over 108}    J_{p}  + \frac{1}{12} I_{p}  &  {1\over 18}  J_{p,q} +\frac{1}{12}\Gamma^T  \\
\cline{1-3}
{1\over 108} J_{p}  +\frac{1}{12} I_{p}  & {5\over 54}  J_p  + \frac{1}{6}I_p & {1\over 18}   J_{p,q} +\frac{1}{12}\Gamma^T   \\
\cline{1-3}
{1\over 18}  J_{q,p} +\frac{1}{12}\Gamma & {1\over 18} J_{q,p} +\frac{1}{12}\Gamma 			 & \frac{1}{12}I_q  + \frac{1}{12}\Gamma \Gamma^T + \frac{1}{12}J_{ q  }	  \\
\end{block}
\end{blockarray} \ .
\end{gather*}
\endgroup
Setting $B_1= 6\widetilde B_1$, we obtain $B\succeq 0 \Longleftrightarrow \widetilde B_1\succeq 0 \Longleftrightarrow B_1\succeq 0$, where
\begingroup
\fontsize{10pt}{18pt}\selectfont
\begin{gather*}
B_1 =
\begin{blockarray}{ccc}
\begin{block}{(c|c|c)}
{5\over 9}  J_{p} +  I_{p}    & {1\over 18}    J_{p}  + \frac{1}{2} I_{p}  &  {1\over 3}  J_{p,q} +\frac{1}{2}\Gamma^T  \\
\cline{1-3}
{1\over 18} J_{p}  +\frac{1}{2} I_{p}  & {5\over 9}  J_p  + I_p & {1\over 3}   J_{p,q} +\frac{1}{2}\Gamma^T   \\
\cline{1-3}
{1\over 3}  J_{q,p} +\frac{1}{2}\Gamma & {1\over 3} J_{q,p} +\frac{1}{2}\Gamma 			 & \frac{1}{2}I_q  + \frac{1}{2}\Gamma \Gamma^T + \frac{1}{2}J_{ q  }	  \\
\end{block}
\end{blockarray}\ .
\end{gather*}
\endgroup

\medskip\noindent
{\bf Step 2:} We now take the Schur complement with respect to the upper left corner of $B_1$ (indexed by $I_1$), where we use Lemma \ref{leminverse} to invert it:
$$
(I_{p} + 5/9 J_{p} )^{-1} =  I_p - 5/(5p + 9) J_p.
$$
After taking this Schur complement  the resulting matrix  $\widetilde B_2$   reads:
 \begingroup
\fontsize{10pt}{18pt}\selectfont
$$
\widetilde B_2=\begin{blockarray}{cc}
\begin{block}{(c|c)}
{5\over 9}  J_p  + I_q & {1\over 3}   J_{p,q} +\frac{1}{2}\Gamma^T   \\
\cline{1-2}
 {1\over 3} J_{q,p} +\frac{1}{2}\Gamma 			 & \frac{1}{2}I_q  + \frac{1}{2}\Gamma \Gamma^T + \frac{1}{2}J_{ q }	  \\
\end{block}
\end{blockarray}
\vspace*{-0.5cm}$$
$$
-
\begin{pmatrix}
{1\over 18}    J_{p}  + \frac{1}{2} I_{p}\\
{1\over 3}  J_{q,p} +\frac{1}{2}\Gamma
\end{pmatrix}
\Big( I_p - {5\over (5p + 9)} J_p\Big)
\begin{pmatrix}
{ 1\over 18}    J_{p}  + \frac{1}{2} I_{p}  & { 1\over 3}  J_{p,q} +\frac{1}{2}\Gamma^T\\
\end{pmatrix}
$$
\endgroup	
\begingroup
\fontsize{10pt}{18pt}\selectfont	
$$
= \begin{blockarray}{cc}
\begin{block}{(c|c)}
\frac{3}{4}I_p + \frac{11p + 23}{4(5p+9)} J_p & \frac{1}{4}\Gamma^T + \frac{3p + 7}{2(5p+9)}J_{p,q}\\
\cline{1-2}
 \frac{1}{4}\Gamma + \frac{3p + 7}{2(5p+9)}J_{q,p}    & \frac{1}{2}I_q + \frac{1}{4}\Gamma \Gamma^T + \frac{3p + 7}{2(5p+9)} J_q\\
\end{block}
\end{blockarray}\ .
$$
\endgroup	
Setting $B_2=4\widetilde B_2$ we obtain $B\succeq 0\Longleftrightarrow B_1\succeq 0 \Longleftrightarrow B_2\succeq 0$, where
\begingroup
\fontsize{10pt}{18pt}\selectfont
$$
B_2 =
\begin{blockarray}{cc}
\begin{block}{(c|c)}
3I_p + \frac{11p + 23}{5p+9} J_p & \Gamma^T + \frac{2(3p + 7)}{5p+9}J_{p,q}\\
\cline{1-2}
\Gamma + \frac{2(3p + 7)}{5p+9}J_{q,p}    & 2I_q + \Gamma \Gamma^T + \frac{2(3p + 7)}{5p+9} J_q\\
\end{block}
\end{blockarray}~.
$$
\endgroup		

\medskip\noindent
{\bf Step 3:}
Inverting the top left block of $B_2$ via Lemma \ref{leminverse} gives
$$
\Big(3I_p + \frac{11p + 23}{5p+9} J_p \Big)^{-1} = \frac{1}{3} I_p - \frac{(11p + 23)}{3(11p^2 + 38p + 27)} J_p.
$$
Taking the third and final Schur complement with respect to this block in $B_2$ we get the matrix
\begingroup
\fontsize{10pt}{18pt}\selectfont$$
B_3 :=  2I_q + \Gamma \Gamma^T + \frac{2(3p + 7)}{5p+9} J_q $$
$$
- \Big(\Gamma^T + \frac{2(3p + 7)}{5p+9}J_{q,p} \Big) \Big( \frac{1}{3} I_p - \frac{(11p + 23)}{3(11p^2 + 38p + 27)} J_p\Big) \Big(\Gamma^T + \frac{2(3p + 7)}{5p+9}J_{p,q}\Big)
$$
\endgroup
$$
 = 2I_q +  \frac{2}{3}\Gamma\Gamma^T + \frac{2(9p + 25)}{3(11p + 27)}J_q.
$$
It is now clear that $B_3\succeq 0$. In turn, this implies that $B_2\succeq 0$ and thus $B\succeq 0$, which concludes the proof of Lemma \ref{lemmainf3}.

\subsection{Some numerical justification for convexity of $f_d$}\label{secnumeric}

We have carried out some numerical experiments for a range of values of $d,L,n$ and verified that the matrices $Q_\gamma$ are positive semidefinite for all $\gamma\in \oN^n_{d-2}$ in these cases. Hence for these values the polynomial $f_d$ is convex and  Conjecture \ref{conjfd} holds.
Recall from Lemma \ref{lembahat} that the matrix $Q_\gamma$ can be decomposed as
$$Q_\gamma=M_\gamma\circ\Big( \underbrace{  \sum_{k\in [m]: \gamma_k\ge 1} Q_{\gamma-u_k}}_{=:B_\gamma} +R_\gamma   \Big)= M_\gamma\circ (B_\gamma +R_\gamma).$$

By the results in Section \ref{secprooftheopd} we already know that the matrix $M_\gamma$ is positive semidefinite. Hence it now suffices to show that the matrix $B_\gamma+R_\gamma$ is positive semidefinite for each $\gamma\in \oN^n_{d-2}$.  In Tables \ref{tabled3L2}-\ref{tabled4L4} in Appendix A   we provide information about the minimum eigenvalues of the matrices $Q_\gamma$, $B_\gamma$ and $R_\gamma$
 for different values of $n$, $d$ and $L$. In each case we consider the possible different cases for selecting $\gamma\in \oN^n_{d-2}$ up to symmetry; the different instances of $\gamma$ are indicated in the column labeled $\gamma$.  For instance, for $d=3,L=2$ there is only one possibility, say $\gamma=u_1$ corresponding to edge $e_1=\{1,2\}$ (see Table \ref{tabled3L2}). For $d=4,L=2$ there are three possibilities:
 $\gamma=2e_1$ with $e_1=\{1,2\}$, $\gamma=u_1+u_2$ with $e_1=\{1,2\}$ and $e_2=\{1,3\}$, and $\gamma=u_1+u_2$ with $e_1=\{1,2\}$ and $e_2=\{3,4\}$ (see Table \ref{tabled4L2}).

  In all cases we  find that $Q_\gamma$ is positive semidefinite (in fact, positive definite). As mentioned above for the case $d=3$,  we see that in general this cannot be deduced by considering its summands separately, since $R_\gamma$ has a negative smallest eigenvalue and 
 $\lambda_{min}(B_\gamma) + \lambda_{\min}(R_\gamma)<0$ from a certain $n$  (which depends on $d$ and $L$). In addition  we observe that $\lambda_{\min}(B_\gamma)$ stays constant from a certain $n$ while $\lambda_{\min}(R_\gamma)$ keeps decreasing. 
 


\subsection*{Acknowledgments}
The first author  is supported by the Europeans Union's EU Framework Programme for Research and Innovation Horizon 2020 under the Marie Sk{ł}odowska-Curie Actions Grant Agreement No 764759 (MINOA) and the third author is supported by the Europeans Union's EU Framework Programme for Research and Innovation Horizon 2020 under the Marie Sk{ł}odowska-Curie Actions Grant Agreement No 813211 (POEMA).

\newpage
\appendix
\section{Numerical results for the polynomials $f_d$}

We group here Tables \ref{tabled3L2}-\ref{tabled4L4} which show the eigenvalues of the matrices $Q_\gamma$, $B_\gamma$ and $R_\gamma$ for small  values of $n,d,L$.

\ignore{
\begin{table}
\begin{center}
	\begin{tabular}{||c c c c c c c||}
		\hline
		$d$ & $L$ & $n$ & $\gamma$ & $\lambda_{min}(Q_{\gamma})$ & $\lambda_{min}(B_{\gamma})$ & $\lambda_{min}(R_{\gamma})$  \\ [0.5ex]
		\hline\hline
        2 & 2 & 3 & [] & 0.1667 & 0.0 & -0.0\\
2 & 2 & 4 & [] & 0.0833 & 0.0 & -0.0\\
2 & 2 & 5 & [] & 0.0833 & 0.0 & 0.0\\
2 & 2 & 6 & [] & 0.0833 & 0.0 & -0.0\\
\hline
	\end{tabular}
	\caption{Case  $d=2,L=2$\label{tabled2L2}}
\end{center}	
	\end{table}
}

\begin{table}[H]
\begin{center}
	\begin{tabular}{||c c c c c c c||}
		\hline
		$d$ & $L$ & $n$ & $\gamma$ & $\lambda_{min}(Q_{\gamma})$ & $\lambda_{min}(B_{\gamma})$ & $\lambda_{min}(R_{\gamma})$  \\ [0.5ex]
		\hline\hline
        3 & 2 & 3 & [[1, 2]] & 0.0556 & 0.1667 & -0.0236\\
3 & 2 & 4 & [[1, 2]] & 0.0347 & 0.0833 & -0.0478\\
3 & 2 & 5 & [[1, 2]] & 0.0347 & 0.0833 & -0.0729\\
3 & 2 & 6 & [[1, 2]] & 0.0347 & 0.0833 & -0.0987\\
3 & 2 & 7 & [[1, 2]] & 0.0347 & 0.0833 & -0.1249\\
3 & 2 & 8 & [[1, 2]] & 0.0347 & 0.0833 & -0.1514\\
\hline
	\end{tabular}
	\caption{Case  $d=3,L=2$\label{tabled3L2}}
\end{center}	
	\end{table}

\begin{table}[H]
\begin{center}
	\begin{tabular}{||c c c c c c c||}
		\hline
		$d$ & $L$ & $n$ & $\gamma$ & $\lambda_{min}(Q_{\gamma})$ & $\lambda_{min}(B_{\gamma})$ & $\lambda_{min}(R_{\gamma})$  \\ [0.5ex]
		\hline\hline
4 & 2 & 3 & [[1, 2], [1, 2]] & 0.0185 & 0.0556 & -0.0415\\
4 & 2 & 4 & [[1, 2], [1, 2]] & 0.0133 & 0.0347 & -0.0805\\
4 & 2 & 5 & [[1, 2], [1, 2]] & 0.0133 & 0.0347 & -0.1189\\
4 & 2 & 6 & [[1, 2], [1, 2]] & 0.0133 & 0.0347 & -0.1572\\
\hline
4 & 2 & 3 & [[1, 3], [1, 2]] & 0.0593 & 0.1778 & -0.0028\\
4 & 2 & 4 & [[1, 3], [1, 2]] & 0.0238 & 0.0802 & -0.0478\\
4 & 2 & 5 & [[1, 3], [1, 2]] & 0.0214 & 0.0743 & -0.092\\
4 & 2 & 6 & [[1, 3], [1, 2]] & 0.0214 & 0.0741 & -0.1359\\
4 & 2 & 7 & [[1, 3], [1, 2]] & 0.0214 & 0.074 & -0.1798\\
\hline
4 & 2 & 4 & [[3, 4], [1, 2]] & 0.0174 & 0.0694 & -0.0012\\
4 & 2 & 5 & [[3, 4], [1, 2]] & 0.0174 & 0.0694 & -0.029\\
4 & 2 & 6 & [[3, 4], [1, 2]] & 0.0174 & 0.0694 & -0.0565\\
4 & 2 & 7 & [[3, 4], [1, 2]] & 0.0174 & 0.0694 & -0.084\\
4 & 2 & 8 & [[3, 4], [1, 2]] & 0.0174 & 0.0694 & -0.1115\\
4 & 2 & 9 & [[3, 4], [1, 2]] & 0.0174 & 0.0694 & -0.139\\
\hline
	\end{tabular}
	\caption{Case  $d=4,L=2$\label{tabled4L2}}
\end{center}
\end{table}

\begin{table}
\begin{center}
	\begin{tabular}{||c c c c c c c||}
		\hline
		$d$ & $L$ & $n$ & $\gamma$ & $\lambda_{min}(Q_{\gamma})$ & $\lambda_{min}(B_{\gamma})$ & $\lambda_{min}(R_{\gamma})$  \\ [0.5ex]
		\hline\hline
5 & 2 & 3 & [[1, 2], [1, 2], [1, 2]] & 0.0062 & 0.0185 & -0.0425\\
5 & 2 & 4 & [[1, 2], [1, 2], [1, 2]] & 0.0049 & 0.0133 & -0.0804\\
5 & 2 & 5 & [[1, 2], [1, 2], [1, 2]] & 0.0049 & 0.0133 & -0.1163\\
\hline
5 & 2 & 3 & [[1, 3], [1, 2], [1, 2]] & 0.0298 & 0.0894 & -0.0062\\
5 & 2 & 4 & [[1, 3], [1, 2], [1, 2]] & 0.0111 & 0.0396 & -0.0605\\
5 & 2 & 5 & [[1, 3], [1, 2], [1, 2]] & 0.0098 & 0.0358 & -0.112\\
5 & 2 & 6 & [[1, 3], [1, 2], [1, 2]] & 0.0098 & 0.0358 & -0.162\\
\hline
5 & 2 & 4 & [[3, 4], [1, 2], [1, 2]] & 0.0077 & 0.0307 & -0.0085\\
5 & 2 & 5 & [[3, 4], [1, 2], [1, 2]] & 0.0072 & 0.0307 & -0.038\\
5 & 2 & 6 & [[3, 4], [1, 2], [1, 2]] & 0.0067 & 0.0307 & -0.0667\\
5 & 2 & 7 & [[3, 4], [1, 2], [1, 2]] & 0.0067 & 0.0307 & -0.0948\\
\hline
5 & 2 & 4 & [[1, 4], [1, 3], [1, 2]] & 0.0263 & 0.1052 & -0.009\\
5 & 2 & 5 & [[1, 4], [1, 3], [1, 2]] & 0.0162 & 0.0716 & -0.0681\\
5 & 2 & 6 & [[1, 4], [1, 3], [1, 2]] & 0.0151 & 0.0676 & -0.1255\\
5 & 2 & 7 & [[1, 4], [1, 3], [1, 2]] & 0.015 & 0.0675 & -0.1819\\
5 & 2 & 8 & [[1, 4], [1, 3], [1, 2]] & 0.015 & 0.0675 & -0.2374\\
\hline
5 & 2 & 4 & [[2, 4], [1, 3], [1, 2]] & 0.0188 & 0.0753 & -0.0063\\
5 & 2 & 5 & [[2, 4], [1, 3], [1, 2]] & 0.0151 & 0.0678 & -0.0613\\
5 & 2 & 6 & [[2, 4], [1, 3], [1, 2]] & 0.0139 & 0.0635 & -0.1147\\
5 & 2 & 7 & [[2, 4], [1, 3], [1, 2]] & 0.0139 & 0.0635 & -0.167\\
5 & 2 & 8 & [[2, 4], [1, 3], [1, 2]] & 0.0139 & 0.0635 & -0.2186\\
\hline
5 & 2 & 5 & [[2, 3], [1, 5], [1, 4]] & 0.0114 & 0.0571 & -0.0053\\
5 & 2 & 6 & [[2, 3], [1, 5], [1, 4]] & 0.0113 & 0.0569 & -0.0395\\
5 & 2 & 7 & [[2, 3], [1, 5], [1, 4]] & 0.0107 & 0.0569 & -0.0731\\
5 & 2 & 8 & [[2, 3], [1, 5], [1, 4]] & 0.0107 & 0.0569 & -0.1062\\
5 & 2 & 9 & [[2, 3], [1, 5], [1, 4]] & 0.0107 & 0.0569 & -0.1391\\
\hline
5 & 2 & 3 & [[2, 3], [1, 3], [1, 2]] & 0.0926 & 0.2778 & -0.0\\
5 & 2 & 4 & [[2, 3], [1, 3], [1, 2]] & 0.0237 & 0.085 & -0.0967\\
5 & 2 & 5 & [[2, 3], [1, 3], [1, 2]] & 0.0212 & 0.0764 & -0.1882\\
5 & 2 & 6 & [[2, 3], [1, 3], [1, 2]] & 0.0212 & 0.0764 & -0.2766\\
\hline
5 & 2 & 6 & [[5, 6], [3, 4], [1, 2]] & 0.0087 & 0.0521 & -0.0011\\
5 & 2 & 7 & [[5, 6], [3, 4], [1, 2]] & 0.0087 & 0.0521 & -0.0233\\
5 & 2 & 8 & [[5, 6], [3, 4], [1, 2]] & 0.0087 & 0.0521 & -0.0452\\
5 & 2 & 9 & [[5, 6], [3, 4], [1, 2]] & 0.0087 & 0.0521 & -0.067\\
5 & 2 & 10 & [[5, 6], [3, 4], [1, 2]] & 0.0087 & 0.0521 & -0.0885\\
5 & 2 & 11 & [[5, 6], [3, 4], [1, 2]] & 0.0087 & 0.0521 & -0.11\\
\hline
	\end{tabular}
	\caption{Case  $d=5,L=2$ \label{tabled5L2}}
\end{center}
\end{table}

\begin{table}
\begin{center}
	\begin{tabular}{||c c c c c c c||}
		\hline
		$d$ & $L$ & $n$ & $\gamma$ & $\lambda_{min}(Q_{\gamma})$ & $\lambda_{min}(B_{\gamma})$ & $\lambda_{min}(R_{\gamma})$  \\ [0.5ex]
		\hline\hline
6 & 2 & 3 & [[1, 2], [1, 2], [1, 2], [1, 2]] & 0.0021 & 0.0062 & -0.0349\\
6 & 2 & 4 & [[1, 2], [1, 2], [1, 2], [1, 2]] & 0.0017 & 0.0049 & -0.0652\\
6 & 2 & 5 & [[1, 2], [1, 2], [1, 2], [1, 2]] & 0.0017 & 0.0049 & -0.0931\\
\hline
6 & 2 & 3 & [[1, 3], [1, 2], [1, 2], [1, 2]] & 0.0124 & 0.0371 & -0.0094\\
6 & 2 & 4 & [[1, 3], [1, 2], [1, 2], [1, 2]] & 0.0044 & 0.0165 & -0.0579\\
6 & 2 & 5 & [[1, 3], [1, 2], [1, 2], [1, 2]] & 0.004 & 0.0148 & -0.1029\\
6 & 2 & 6 & [[1, 3], [1, 2], [1, 2], [1, 2]] & 0.004 & 0.0148 & -0.1457\\
\hline
6 & 2 & 3 & [[1, 3], [1, 3], [1, 2], [1, 2]] & 0.0261 & 0.0785 & -0.0016\\
6 & 2 & 4 & [[1, 3], [1, 3], [1, 2], [1, 2]] & 0.0064 & 0.0237 & -0.0626\\
6 & 2 & 5 & [[1, 3], [1, 3], [1, 2], [1, 2]] & 0.0057 & 0.0211 & -0.1193\\
6 & 2 & 6 & [[1, 3], [1, 3], [1, 2], [1, 2]] & 0.0057 & 0.0211 & -0.1732\\
\hline
6 & 2 & 4 & [[3, 4], [1, 2], [1, 2], [1, 2]] & 0.0031 & 0.0125 & -0.0141\\
6 & 2 & 5 & [[3, 4], [1, 2], [1, 2], [1, 2]] & 0.0026 & 0.0124 & -0.0384\\
6 & 2 & 6 & [[3, 4], [1, 2], [1, 2], [1, 2]] & 0.0024 & 0.0115 & -0.0616\\
\hline
6 & 2 & 4 & [[3, 4], [3, 4], [1, 2], [1, 2]] & 0.0038 & 0.0153 & -0.0007\\
6 & 2 & 5 & [[3, 4], [3, 4], [1, 2], [1, 2]] & 0.0036 & 0.0153 & -0.0255\\
6 & 2 & 6 & [[3, 4], [3, 4], [1, 2], [1, 2]] & 0.0033 & 0.0153 & -0.0492\\
6 & 2 & 7 & [[3, 4], [3, 4], [1, 2], [1, 2]] & 0.0033 & 0.0153 & -0.0721\\
\hline
6 & 2 & 4 & [[1, 4], [1, 3], [1, 2], [1, 2]] & 0.0147 & 0.0589 & -0.0142\\
6 & 2 & 5 & [[1, 4], [1, 3], [1, 2], [1, 2]] & 0.0084 & 0.0386 & -0.0771\\
6 & 2 & 6 & [[1, 4], [1, 3], [1, 2], [1, 2]] & 0.0078 & 0.036 & -0.1372\\
6 & 2 & 7 & [[1, 4], [1, 3], [1, 2], [1, 2]] & 0.0078 & 0.036 & -0.1952\\
\hline
6 & 2 & 4 & [[2, 4], [1, 3], [1, 2], [1, 2]] & 0.0129 & 0.0514 & -0.0151\\
6 & 2 & 5 & [[2, 4], [1, 3], [1, 2], [1, 2]] & 0.0079 & 0.037 & -0.0766\\
6 & 2 & 6 & [[2, 4], [1, 3], [1, 2], [1, 2]] & 0.0073 & 0.0344 & -0.1352\\
6 & 2 & 7 & [[2, 4], [1, 3], [1, 2], [1, 2]] & 0.0073 & 0.0344 & -0.1919\\
\hline
6 & 2 & 4 & [[2, 4], [1, 3], [1, 3], [1, 2]] & 0.0102 & 0.0407 & -0.0089\\
6 & 2 & 5 & [[2, 4], [1, 3], [1, 3], [1, 2]] & 0.0074 & 0.0343 & -0.064\\
6 & 2 & 6 & [[2, 4], [1, 3], [1, 3], [1, 2]] & 0.0068 & 0.0318 & -0.1167\\
6 & 2 & 7 & [[2, 4], [1, 3], [1, 3], [1, 2]] & 0.0068 & 0.0318 & -0.1675\\
\hline
6 & 2 & 5 & [[2, 3], [1, 5], [1, 4], [1, 4]] & 0.0059 & 0.0294 & -0.0148\\
6 & 2 & 6 & [[2, 3], [1, 5], [1, 4], [1, 4]] & 0.0052 & 0.0293 & -0.0485\\
6 & 2 & 7 & [[2, 3], [1, 5], [1, 4], [1, 4]] & 0.0049 & 0.0278 & -0.0813\\
6 & 2 & 8 & [[2, 3], [1, 5], [1, 4], [1, 4]] & 0.0049 & 0.0278 & -0.1132\\
\hline
6 & 2 & 5 & [[2, 3], [2, 3], [1, 5], [1, 4]] & 0.0053 & 0.0266 & -0.0055\\
6 & 2 & 6 & [[2, 3], [2, 3], [1, 5], [1, 4]] & 0.0047 & 0.0259 & -0.0344\\
6 & 2 & 7 & [[2, 3], [2, 3], [1, 5], [1, 4]] & 0.0044 & 0.0249 & -0.0623\\
6 & 2 & 8 & [[2, 3], [2, 3], [1, 5], [1, 4]] & 0.0044 & 0.0249 & -0.0897\\
\hline
\hline
\end{tabular}
\end{center}
\end{table}

\begin{table}
\begin{center}
	\begin{tabular}{||c c c c c c c||}
\hline
6 & 2 & 3 & [[2, 3], [1, 3], [1, 2], [1, 2]] & 0.0525 & 0.1574 & -0.0017\\
6 & 2 & 4 & [[2, 3], [1, 3], [1, 2], [1, 2]] & 0.0125 & 0.0467 & -0.1176\\
6 & 2 & 5 & [[2, 3], [1, 3], [1, 2], [1, 2]] & 0.0112 & 0.0416 & -0.2252\\
6 & 2 & 6 & [[2, 3], [1, 3], [1, 2], [1, 2]] & 0.0112 & 0.0416 & -0.3274\\
\hline
6 & 2 & 6 & [[5, 6], [3, 4], [1, 2], [1, 2]] & 0.0038 & 0.023 & -0.0086\\
6 & 2 & 7 & [[5, 6], [3, 4], [1, 2], [1, 2]] & 0.0035 & 0.0229 & -0.0278\\
6 & 2 & 8 & [[5, 6], [3, 4], [1, 2], [1, 2]] & 0.0033 & 0.022 & -0.0466\\
6 & 2 & 9 & [[5, 6], [3, 4], [1, 2], [1, 2]] & 0.0033 & 0.022 & -0.065\\
		\hline
		6 & 2 & 5 & [[1, 5], [1, 4], [1, 3], [1, 2]] & 0.0208 & 0.104 & -0.0167\\
6 & 2 & 6 & [[1, 5], [1, 4], [1, 3], [1, 2]] & 0.0121 & 0.066 & -0.0852\\
6 & 2 & 7 & [[1, 5], [1, 4], [1, 3], [1, 2]] & 0.0115 & 0.063 & -0.1515\\
6 & 2 & 8 & [[1, 5], [1, 4], [1, 3], [1, 2]] & 0.0115 & 0.0629 & -0.2164\\
\hline
6 & 2 & 5 & [[2, 3], [1, 5], [1, 4], [1, 2]] & 0.0141 & 0.0706 & -0.0132\\
6 & 2 & 6 & [[2, 3], [1, 5], [1, 4], [1, 2]] & 0.0107 & 0.0592 & -0.0743\\
6 & 2 & 7 & [[2, 3], [1, 5], [1, 4], [1, 2]] & 0.0101 & 0.0562 & -0.1336\\
6 & 2 & 8 & [[2, 3], [1, 5], [1, 4], [1, 2]] & 0.0101 & 0.0562 & -0.1915\\
\hline
6 & 2 & 6 & [[2, 3], [1, 6], [1, 5], [1, 4]] & 0.0084 & 0.0505 & -0.0119\\
6 & 2 & 7 & [[2, 3], [1, 6], [1, 5], [1, 4]] & 0.0079 & 0.0503 & -0.0503\\
6 & 2 & 8 & [[2, 3], [1, 6], [1, 5], [1, 4]] & 0.0075 & 0.049 & -0.0879\\
6 & 2 & 9 & [[2, 3], [1, 6], [1, 5], [1, 4]] & 0.0075 & 0.0489 & -0.1248\\
\hline
6 & 2 & 4 & [[2, 3], [1, 4], [1, 3], [1, 2]] & 0.0246 & 0.0985 & -0.0122\\
6 & 2 & 5 & [[2, 3], [1, 4], [1, 3], [1, 2]] & 0.0162 & 0.0746 & -0.1185\\
6 & 2 & 6 & [[2, 3], [1, 4], [1, 3], [1, 2]] & 0.0151 & 0.0695 & -0.2198\\
6 & 2 & 7 & [[2, 3], [1, 4], [1, 3], [1, 2]] & 0.0151 & 0.0695 & -0.3177\\
\hline
6 & 2 & 4 & [[3, 4], [2, 4], [1, 3], [1, 2]] & 0.0204 & 0.0815 & -0.0003\\
6 & 2 & 5 & [[3, 4], [2, 4], [1, 3], [1, 2]] & 0.014 & 0.0644 & -0.0946\\
6 & 2 & 6 & [[3, 4], [2, 4], [1, 3], [1, 2]] & 0.0129 & 0.0594 & -0.1846\\
6 & 2 & 7 & [[3, 4], [2, 4], [1, 3], [1, 2]] & 0.0129 & 0.0594 & -0.2716\\
\hline
6 & 2 & 6 & [[2, 6], [2, 3], [1, 5], [1, 4]] & 0.0077 & 0.046 & -0.005\\
6 & 2 & 7 & [[2, 6], [2, 3], [1, 5], [1, 4]] & 0.0074 & 0.0456 & -0.0392\\
6 & 2 & 8 & [[2, 6], [2, 3], [1, 5], [1, 4]] & 0.0071 & 0.0456 & -0.0727\\
6 & 2 & 9 & [[2, 6], [2, 3], [1, 5], [1, 4]] & 0.0071 & 0.0456 & -0.1055\\
\hline
6 & 2 & 5 & [[2, 5], [2, 3], [1, 4], [1, 3]] & 0.0121 & 0.0603 & -0.0084\\
6 & 2 & 6 & [[2, 5], [2, 3], [1, 4], [1, 3]] & 0.01 & 0.055 & -0.0643\\
6 & 2 & 7 & [[2, 5], [2, 3], [1, 4], [1, 3]] & 0.0094 & 0.0523 & -0.1184\\
6 & 2 & 8 & [[2, 5], [2, 3], [1, 4], [1, 3]] & 0.0094 & 0.0523 & -0.1712\\
\hline
6 & 2 & 6 & [[5, 6], [2, 4], [1, 3], [1, 2]] & 0.0077 & 0.0461 & -0.0096\\
6 & 2 & 7 & [[5, 6], [2, 4], [1, 3], [1, 2]] & 0.0073 & 0.0461 & -0.0451\\
6 & 2 & 8 & [[5, 6], [2, 4], [1, 3], [1, 2]] & 0.007 & 0.0457 & -0.0799\\
6 & 2 & 9 & [[5, 6], [2, 4], [1, 3], [1, 2]] & 0.007 & 0.0457 & -0.114\\
\hline
	\end{tabular}
\end{center}
\end{table}

\begin{table}
\begin{center}
	\begin{tabular}{||c c c c c c c||}
\hline
6 & 2 & 7 & [[4, 5], [2, 3], [1, 7], [1, 6]] & 0.0057 & 0.0399 & -0.0047\\
6 & 2 & 8 & [[4, 5], [2, 3], [1, 7], [1, 6]] & 0.0056 & 0.0399 & -0.0276\\
6 & 2 & 9 & [[4, 5], [2, 3], [1, 7], [1, 6]] & 0.0053 & 0.0398 & -0.0501\\
6 & 2 & 10 & [[4, 5], [2, 3], [1, 7], [1, 6]] & 0.0053 & 0.0398 & -0.0723\\
\hline
6 & 2 & 5 & [[4, 5], [2, 3], [1, 3], [1, 2]] & 0.0121 & 0.0606 & -0.0192\\
6 & 2 & 6 & [[4, 5], [2, 3], [1, 3], [1, 2]] & 0.0112 & 0.0606 & -0.0784\\
6 & 2 & 7 & [[4, 5], [2, 3], [1, 3], [1, 2]] & 0.0106 & 0.0594 & -0.1359\\
6 & 2 & 8 & [[4, 5], [2, 3], [1, 3], [1, 2]] & 0.0106 & 0.0594 & -0.1919\\
\hline
6 & 2 & 8 & [[7, 8], [5, 6], [3, 4], [1, 2]] & 0.0043 & 0.0347 & -0.0008\\
6 & 2 & 9 & [[7, 8], [5, 6], [3, 4], [1, 2]] & 0.0043 & 0.0347 & -0.0161\\
6 & 2 & 10 & [[7, 8], [5, 6], [3, 4], [1, 2]] & 0.0043 & 0.0347 & -0.0312\\
6 & 2 & 11 & [[7, 8], [5, 6], [3, 4], [1, 2]] & 0.0043 & 0.0347 & -0.0462\\
\hline
	\end{tabular}
\end{center}
\caption{Case  $d=6, L=2$ \label{tabled6L2}}
\end{table}

\begin{table}
\begin{center}
	\begin{tabular}{||c c c c c c c||}
		\hline
		$d$ & $L$ & $n$ & $\gamma$ & $\lambda_{min}(Q_{\gamma})$ & $\lambda_{min}(B_{\gamma})$ & $\lambda_{min}(R_{\gamma})$  \\ [0.5ex]
		\hline\hline
2 & 3 & 3 &  & 0.2222 & 0.0 & 0.6667\\
2 & 3 & 4 &  & 0.0556 & 0.0 & -0.0\\
2 & 3 & 5 &  & 0.0222 & 0.0 & -0.0\\
2 & 3 & 6 &  & 0.0111 & 0.0 & -0.0\\
\hline
	\end{tabular}
\end{center}
\caption{Case $d=2,L=3$ \label{tabled2L3}}
\end{table}

\begin{table}\begin{center}
	\begin{tabular}{||c c c c c c c||}
		\hline
		$d$ & $L$ & $n$ & $\gamma$ & $\lambda_{min}(Q_{\gamma})$ & $\lambda_{min}(B_{\gamma})$ & $\lambda_{min}(R_{\gamma})$  \\ [0.5ex]
		\hline\hline
3 & 3 & 3 & [[1, 2, 3]] & 0.2222 & 0.2222 & 0.4444\\
3 & 3 & 4 & [[1, 2, 3]] & 0.0139 & 0.0556 & -0.0064\\
3 & 3 & 5 & [[1, 2, 3]] & 0.0053 & 0.0222 & -0.0191\\
3 & 3 & 6 & [[1, 2, 3]] & 0.0033 & 0.0111 & -0.0382\\
3 & 3 & 7 & [[1, 2, 3]] & 0.0032 & 0.0111 & -0.064\\
3 & 3 & 8 & [[1, 2, 3]] & 0.0032 & 0.0111 & -0.0963\\
\hline
	\end{tabular}
	\caption{Case $d=3,L=3$ \label{tabled3L3}}
\end{center}
\end{table}

\begin{table}
\begin{center}
	\begin{tabular}{||c c c c c c c||}
		\hline
		$d$ & $L$ & $n$ & $\gamma$ & $\lambda_{min}(Q_{\gamma})$ & $\lambda_{min}(B_{\gamma})$ & $\lambda_{min}(R_{\gamma})$  \\ [0.5ex]
		\hline\hline
4 & 3 & 3 & [[1, 2, 3], [1, 2, 3]] & 0.1481 & 0.2222 & 0.2222\\
4 & 3 & 4 & [[1, 2, 3], [1, 2, 3]] & 0.0035 & 0.0139 & -0.0081\\
4 & 3 & 5 & [[1, 2, 3], [1, 2, 3]] & 0.0012 & 0.0053 & -0.0216\\
4 & 3 & 6 & [[1, 2, 3], [1, 2, 3]] & 0.0009 & 0.0033 & -0.0409\\
4 & 3 & 7 & [[1, 2, 3], [1, 2, 3]] & 0.0008 & 0.0032 & -0.0659\\
\hline
4 & 3 & 4 & [[1, 2, 3], [1, 2, 4]] & 0.0069 & 0.0278 & -0.0007\\
4 & 3 & 5 & [[1, 2, 3], [1, 2, 4]] & 0.0023 & 0.0114 & -0.0176\\
4 & 3 & 6 & [[1, 2, 3], [1, 2, 4]] & 0.0015 & 0.0067 & -0.0424\\
4 & 3 & 7 & [[1, 2, 3], [1, 2, 4]] & 0.0014 & 0.0065 & -0.0751\\
\hline
4 & 3 & 5 & [[1, 2, 5], [1, 3, 4]] & 0.0025 & 0.0126 & -0.001\\
4 & 3 & 6 & [[1, 2, 5], [1, 3, 4]] & 0.0013 & 0.0066 & -0.0169\\
4 & 3 & 7 & [[1, 2, 5], [1, 3, 4]] & 0.0012 & 0.0063 & -0.0386\\
4 & 3 & 8 & [[1, 2, 5], [1, 3, 4]] & 0.0012 & 0.0063 & -0.0662\\
\hline
4 & 3 & 6 & [[1, 2, 6], [3, 4, 5]] & 0.0011 & 0.0069 & -0.0004\\
4 & 3 & 7 & [[1, 2, 6], [3, 4, 5]] & 0.001 & 0.0065 & -0.0147\\
4 & 3 & 8 & [[1, 2, 6], [3, 4, 5]] & 0.001 & 0.0063 & -0.0335\\
4 & 3 & 9 & [[1, 2, 6], [3, 4, 5]] & 0.001 & 0.0063 & -0.0566\\
\hline
	\end{tabular}
	\caption{Case $d=4,L=3$\label{tabled4L3}}
\end{center}
\end{table}

\begin{table}
\begin{center}
	\begin{tabular}{||c c c c c c c||}
		\hline
		$d$ & $L$ & $n$ & $\gamma$ & $\lambda_{min}(Q_{\gamma})$ & $\lambda_{min}(B_{\gamma})$ & $\lambda_{min}(R_{\gamma})$  \\ [0.5ex]
		\hline\hline
5 & 3 & 3 & [[1, 2, 3], [1, 2, 3], [1, 2, 3]] & 0.0823 & 0.1481 & 0.0988\\
5 & 3 & 4 & [[1, 2, 3], [1, 2, 3], [1, 2, 3]] & 0.0009 & 0.0035 & -0.0059\\
5 & 3 & 5 & [[1, 2, 3], [1, 2, 3], [1, 2, 3]] & 0.0003 & 0.0012 & -0.0146\\
5 & 3 & 6 & [[1, 2, 3], [1, 2, 3], [1, 2, 3]] & 0.0002 & 0.0009 & -0.0262\\
\hline
5 & 3 & 4 & [[1, 2, 3], [1, 2, 4], [1, 2, 4]] & 0.0026 & 0.0104 & -0.0009\\
5 & 3 & 5 & [[1, 2, 3], [1, 2, 4], [1, 2, 4]] & 0.0007 & 0.0037 & -0.016\\
5 & 3 & 6 & [[1, 2, 3], [1, 2, 4], [1, 2, 4]] & 0.0005 & 0.0024 & -0.0366\\
5 & 3 & 7 & [[1, 2, 3], [1, 2, 4], [1, 2, 4]] & 0.0004 & 0.0022 & -0.0626\\
\hline
5 & 3 & 5 & [[1, 2, 5], [1, 3, 4], [1, 3, 4]] & 0.0008 & 0.0038 & -0.0018\\
5 & 3 & 6 & [[1, 2, 5], [1, 3, 4], [1, 3, 4]] & 0.0004 & 0.0021 & -0.0145\\
5 & 3 & 7 & [[1, 2, 5], [1, 3, 4], [1, 3, 4]] & 0.0004 & 0.002 & -0.0308\\
5 & 3 & 8 & [[1, 2, 5], [1, 3, 4], [1, 3, 4]] & 0.0004 & 0.002 & -0.0508\\
\hline
5 & 3 & 6 & [[1, 2, 6], [3, 4, 5], [3, 4, 5]] & 0.0003 & 0.0021 & -0.0027\\
5 & 3 & 7 & [[1, 2, 6], [3, 4, 5], [3, 4, 5]] & 0.0003 & 0.0019 & -0.0134\\
5 & 3 & 8 & [[1, 2, 6], [3, 4, 5], [3, 4, 5]] & 0.0003 & 0.0018 & -0.0267\\
\hline
5 & 3 & 4 & [[1, 2, 3], [1, 2, 4], [1, 3, 4]] & 0.0067 & 0.027 & -0.0001\\
5 & 3 & 5 & [[1, 2, 3], [1, 2, 4], [1, 3, 4]] & 0.0014 & 0.007 & -0.0287\\
5 & 3 & 6 & [[1, 2, 3], [1, 2, 4], [1, 3, 4]] & 0.001 & 0.0046 & -0.0672\\
5 & 3 & 7 & [[1, 2, 3], [1, 2, 4], [1, 3, 4]] & 0.0009 & 0.0043 & -0.1159\\
\hline
5 & 3 & 5 & [[1, 2, 3], [1, 2, 4], [1, 2, 5]] & 0.0018 & 0.0089 & -0.0032\\
5 & 3 & 6 & [[1, 2, 3], [1, 2, 4], [1, 2, 5]] & 0.0009 & 0.0048 & -0.0289\\
5 & 3 & 7 & [[1, 2, 3], [1, 2, 4], [1, 2, 5]] & 0.0007 & 0.0042 & -0.0622\\
5 & 3 & 8 & [[1, 2, 3], [1, 2, 4], [1, 2, 5]] & 0.0007 & 0.0041 & -0.103\\
\hline
5 & 3 & 5 & [[1, 2, 3], [1, 2, 4], [1, 3, 5]] & 0.0016 & 0.008 & -0.0023\\
5 & 3 & 6 & [[1, 2, 3], [1, 2, 4], [1, 3, 5]] & 0.0008 & 0.0044 & -0.0269\\
5 & 3 & 7 & [[1, 2, 3], [1, 2, 4], [1, 3, 5]] & 0.0007 & 0.0042 & -0.0586\\
5 & 3 & 8 & [[1, 2, 3], [1, 2, 4], [1, 3, 5]] & 0.0007 & 0.004 & -0.0975\\
\hline
5 & 3 & 6 & [[1, 2, 3], [1, 2, 4], [1, 5, 6]] & 0.0007 & 0.0042 & -0.0033\\
5 & 3 & 7 & [[1, 2, 3], [1, 2, 4], [1, 5, 6]] & 0.0006 & 0.0039 & -0.0242\\
5 & 3 & 8 & [[1, 2, 3], [1, 2, 4], [1, 5, 6]] & 0.0006 & 0.0037 & -0.0504\\
5 & 3 & 9 & [[1, 2, 3], [1, 2, 4], [1, 5, 6]] & 0.0006 & 0.0037 & -0.0818\\
\hline
5 & 3 & 5 & [[1, 2, 3], [1, 2, 4], [3, 4, 5]] & 0.0015 & 0.0076 & -0.0012\\
5 & 3 & 6 & [[1, 2, 3], [1, 2, 4], [3, 4, 5]] & 0.0007 & 0.0041 & -0.0245\\
5 & 3 & 7 & [[1, 2, 3], [1, 2, 4], [3, 4, 5]] & 0.0007 & 0.0038 & -0.0545\\
5 & 3 & 8 & [[1, 2, 3], [1, 2, 4], [3, 4, 5]] & 0.0007 & 0.0038 & -0.0913\\
\hline
5 & 3 & 6 & [[1, 2, 3], [1, 2, 4], [3, 5, 6]] & 0.0007 & 0.004 & -0.0027\\
5 & 3 & 7 & [[1, 2, 3], [1, 2, 4], [3, 5, 6]] & 0.0006 & 0.0036 & -0.0228\\
5 & 3 & 8 & [[1, 2, 3], [1, 2, 4], [3, 5, 6]] & 0.0005 & 0.0035 & -0.048\\
5 & 3 & 9 & [[1, 2, 3], [1, 2, 4], [3, 5, 6]] & 0.0005 & 0.0035 & -0.0781\\

\hline
	\end{tabular}
\end{center}
\end{table}

\begin{table}
\begin{center}
	\begin{tabular}{||c c c c c c c||}
		\hline
5 & 3 & 6 & [[1, 2, 4], [1, 3, 5], [2, 3, 6]] & 0.0007 & 0.004 & -0.0019\\
5 & 3 & 7 & [[1, 2, 4], [1, 3, 5], [2, 3, 6]] & 0.0006 & 0.0038 & -0.0216\\
5 & 3 & 8 & [[1, 2, 4], [1, 3, 5], [2, 3, 6]] & 0.0006 & 0.0036 & -0.0462\\
5 & 3 & 9 & [[1, 2, 4], [1, 3, 5], [2, 3, 6]] & 0.0006 & 0.0036 & -0.0758\\
\hline
5 & 3 & 7 & [[1, 2, 5], [1, 3, 4], [1, 6, 7]] & 0.0005 & 0.0035 & -0.003\\
5 & 3 & 8 & [[1, 2, 5], [1, 3, 4], [1, 6, 7]] & 0.0005 & 0.0035 & -0.0207\\
5 & 3 & 9 & [[1, 2, 5], [1, 3, 4], [1, 6, 7]] & 0.0005 & 0.0035 & -0.0423\\
5 & 3 & 10 & [[1, 2, 5], [1, 3, 4], [1, 6, 7]] & 0.0005 & 0.0035 & -0.0678\\
\hline
5 & 3 & 7 & [[1, 2, 7], [1, 3, 4], [2, 5, 6]] & 0.0005 & 0.0034 & -0.0023\\
5 & 3 & 8 & [[1, 2, 7], [1, 3, 4], [2, 5, 6]] & 0.0005 & 0.0034 & -0.0192\\
5 & 3 & 9 & [[1, 2, 7], [1, 3, 4], [2, 5, 6]] & 0.0005 & 0.0034 & -0.0399\\
5 & 3 & 10 & [[1, 2, 7], [1, 3, 4], [2, 5, 6]] & 0.0005 & 0.0034 & -0.0643\\
\hline
5 & 3 & 7 & [[1, 2, 3], [1, 2, 4], [5, 6, 7]] & 0.0005 & 0.0037 & -0.0033\\
5 & 3 & 8 & [[1, 2, 3], [1, 2, 4], [5, 6, 7]] & 0.0005 & 0.0034 & -0.0207\\
5 & 3 & 9 & [[1, 2, 3], [1, 2, 4], [5, 6, 7]] & 0.0004 & 0.0033 & -0.0419\\
5 & 3 & 10 & [[1, 2, 3], [1, 2, 4], [5, 6, 7]] & 0.0004 & 0.0033 & -0.0669\\
\hline
5 & 3 & 8 & [[1, 2, 5], [1, 3, 4], [6, 7, 8]] & 0.0004 & 0.0032 & -0.0018\\
5 & 3 & 9 & [[1, 2, 5], [1, 3, 4], [6, 7, 8]] & 0.0004 & 0.0032 & -0.0163\\
5 & 3 & 10 & [[1, 2, 5], [1, 3, 4], [6, 7, 8]] & 0.0004 & 0.0032 & -0.0336\\
5 & 3 & 11 & [[1, 2, 5], [1, 3, 4], [6, 7, 8]] & 0.0004 & 0.0032 & -0.0539\\
\hline
5 & 3 & 9 & [[1, 2, 6], [3, 4, 5], [7, 8, 9]] & 0.0003 & 0.003 & -0.0005\\
5 & 3 & 10 & [[1, 2, 6], [3, 4, 5], [7, 8, 9]] & 0.0003 & 0.003 & -0.0128\\
5 & 3 & 11 & [[1, 2, 6], [3, 4, 5], [7, 8, 9]] & 0.0003 & 0.003 & -0.0273\\
5 & 3 & 12 & [[1, 2, 6], [3, 4, 5], [7, 8, 9]] & 0.0003 & 0.003 & -0.044\\
\hline
	\end{tabular}
	\caption{Case $d=5,L=3$\label{tabled5L3}}
\end{center}
\end{table}

\ignore{
\begin{table}
\begin{center}
	\begin{tabular}{||c c c c c c c||}
		\hline
		$d$ & $L$ & $n$ & $\gamma$ & $\lambda_{min}(Q_{\gamma})$ & $\lambda_{min}(B_{\gamma})$ & $\lambda_{min}(R_{\gamma})$  \\ [0.5ex]
		\hline\hline
2 & 4 & 4 &  & 0.125 & 0.0 & 0.5\\
2 & 4 & 5 &  & 0.025 & 0.0 & -0.0\\
2 & 4 & 6 &  & 0.0083 & 0.0 & -0.0\\
2 & 4 & 7 &  & 0.0036 & 0.0 & -0.0\\
\hline
	\end{tabular}
\end{center}
\end{table}
}

\begin{table}
\begin{center}
	\begin{tabular}{||c c c c c c c||}
		\hline
		$d$ & $L$ & $n$ & $\gamma$ & $\lambda_{min}(Q_{\gamma})$ & $\lambda_{min}(B_{\gamma})$ & $\lambda_{min}(R_{\gamma})$  \\ [0.5ex]
		\hline\hline
3 & 4 & 4 & [[1, 2, 3, 4]] & 0.0938 & 0.125 & 0.25\\
3 & 4 & 5 & [[1, 2, 3, 4]] & 0.005 & 0.025 & -0.0024\\
3 & 4 & 6 & [[1, 2, 3, 4]] & 0.0014 & 0.0083 & -0.0101\\
3 & 4 & 7 & [[1, 2, 3, 4]] & 0.0007 & 0.0036 & -0.0256\\
3 & 4 & 8 & [[1, 2, 3, 4]] & 0.0004 & 0.0018 & -0.0514\\
\hline
	\end{tabular}
	\caption{Case $d=3,L=4$ \label{tabled3L4}}
\end{center}
\end{table}

\begin{table}
\begin{center}
	\begin{tabular}{||c c c c c c c||}
		\hline
		$d$ & $L$ & $n$ & $\gamma$ & $\lambda_{min}(Q_{\gamma})$ & $\lambda_{min}(B_{\gamma})$ & $\lambda_{min}(R_{\gamma})$  \\ [0.5ex]
		\hline\hline
4 & 4 & 4 & [[1, 2, 3, 4], [1, 2, 3, 4]] & 0.0469 & 0.0938 & 0.0938\\
4 & 4 & 5 & [[1, 2, 3, 4], [1, 2, 3, 4]] & 0.001 & 0.005 & -0.0023\\
4 & 4 & 6 & [[1, 2, 3, 4], [1, 2, 3, 4]] & 0.0002 & 0.0014 & -0.0087\\
4 & 4 & 7 & [[1, 2, 3, 4], [1, 2, 3, 4]] & 0.0001 & 0.0007 & -0.0207\\
4 & 4 & 8 & [[1, 2, 3, 4], [1, 2, 3, 4]] & 0.0001 & 0.0004 & -0.0399\\
\hline
4 & 4 & 5 & [[1, 2, 3, 4], [1, 2, 3, 5]] & 0.002 & 0.01 & -0.0002\\
4 & 4 & 6 & [[1, 2, 3, 4], [1, 2, 3, 5]] & 0.0005 & 0.003 & -0.0081\\
4 & 4 & 7 & [[1, 2, 3, 4], [1, 2, 3, 5]] & 0.0002 & 0.0014 & -0.0242\\
4 & 4 & 8 & [[1, 2, 3, 4], [1, 2, 3, 5]] & 0.0001 & 0.0008 & -0.051\\
\hline
4 & 4 & 6 & [[1, 2, 3, 6], [1, 2, 4, 5]] & 0.0005 & 0.0032 & -0.0006\\
4 & 4 & 7 & [[1, 2, 3, 6], [1, 2, 4, 5]] & 0.0002 & 0.0014 & -0.0109\\
4 & 4 & 8 & [[1, 2, 3, 6], [1, 2, 4, 5]] & 0.0001 & 0.0008 & -0.0292\\
4 & 4 & 9 & [[1, 2, 3, 6], [1, 2, 4, 5]] & 0.0001 & 0.0008 & -0.0575\\
\hline
4 & 4 & 7 & [[1, 2, 3, 7], [1, 4, 5, 6]] & 0.0002 & 0.0016 & -0.0007\\
4 & 4 & 8 & [[1, 2, 3, 7], [1, 4, 5, 6]] & 0.0001 & 0.0008 & -0.0127\\
4 & 4 & 9 & [[1, 2, 3, 7], [1, 4, 5, 6]] & 0.0001 & 0.0008 & -0.0323\\
4 & 4 & 10 & [[1, 2, 3, 7], [1, 4, 5, 6]] & 0.0001 & 0.0008 & -0.0612\\
\hline
4 & 4 & 8 & [[1, 2, 3, 8], [4, 5, 6, 7]] & 0.0001 & 0.0008 & -0.0003\\
4 & 4 & 9 & [[1, 2, 3, 8], [4, 5, 6, 7]] & 0.0001 & 0.0008 & -0.0133\\
4 & 4 & 10 & [[1, 2, 3, 8], [4, 5, 6, 7]] & 0.0001 & 0.0008 & -0.0335\\
\hline
	\end{tabular}
	\caption{Case $d=4,L=4$\label{tabled4L4}}
\end{center}
\end{table}



\begin{thebibliography}{10}

\bibitem{BGSV}
C. Bachoc, D.C. Gijswijt, A. Schrijver, and F. Vallentin.
 Invariant semidefinite programs. 
 In {\em Handbook on Semidefinite, Conic and Polynomial Optimization} (M.F. Anjos, J.B. Lasserre (eds.)), pages 219--269, Springer, 2012.
 
 \bibitem{BEK}
 G. P. Barker, L. Q. Eifler, and T. P. Kezlan. A non-commutative spectral theorem.
 {\em Linear Algebra and Appl.,} 20(2):95--100, 1978.
 
 
 \bibitem{CBS}
E. Cardinaels, S. Borst and J.S.H. van Leeuwaarden.
Redundancy scheduling with locally stable compatibility graphs.
arXiv:2005.14566v1, 2020.

\bibitem{dKPS}
E. de Klerk, D. Pasechnik, and A. Schrijver.
Reductions of symmetric semidefinite programs using the regular $*$-representation.
{\em Mathematical Programming},  109:613--624, 2007.


\bibitem{Delsarte73}
P. Delsarte. 
An Algebraic Approach to the Association Schemes of Coding Theory.
[Phiips Research Reports Supplements (1973) No. 10] 
Philips Research Laboratories.


\bibitem{DelsarteLevenshtein}
P. Delsarte and V. I. Levenshtein. Association schemes and coding theory. {\em IEEE Transactions on  Information Theory},  44(6):2477--2504, 1998. 


\bibitem{GZDHBHSW}
K. Gardner, S. Zbarsky, S. Doroudi, M. Harchol-Balter, E. Hyyti\"a, and A. Shceller-Wolf.
Queueing with redundant requests: Exact analysis.
{\em Queueing Systems,} 83(3--4):227--259, 2016.

\bibitem{GatermanParrilo}
K. Gatermann and P.A. Parrilo. Symmetry groups, semidefinite programs, and sums of squares.
{\em Journal of  Pure and Applied Algebra},  192(1-3):9--128, 2004.

\bibitem{GST}
D. Gijswijt, A. Schrijver, and H. Tanaka.
 New upper bounds for nonbinary codes based on the Terwilliger algebra and semidefinite programming.
 {\em  Journal of Comb. Theory, Series A}, 113(8):1719--1731, 2008.
 
 \bibitem{GMS12}
D. Gijswijt, H.D. Mittelmann,  and A. Schrijver. 
Semidefinite code bounds based on quadruple distances. 
{\em IEEE Transactions on Information Theory},  58:2697--2705, 2012.

\bibitem{Laurent2007}
M. Laurent.  Strengthened semidefinite programming bounds for codes. 
{\em Mathematical Programming,} 109(2-3):239--261, 2007. 

\bibitem{LPS2017}
B. Litjens, S. Polak, and A. Schrijver.
Semidefinite bounds for nonbinary codes based on quadruples.
{\em Codes, Designs and Cryptography}, 84:87--100, 2017.


\bibitem{RSST2018} 
A. Raymond, J. Saunderson, M. Singh, and R. Thomas.  Symmetric sums of squares over k-subset hypercubes.
{\em Mathematical Programming Series A}, 167(2):315--354, 2018.

\bibitem{RSTTuran2018} 
A. Raymond,  M. Singh, and R. Thomas. Symmetry in Tur\'an sums of squares polynomials from flag algebras.
{\em Algebraic Combinatorics}, 1(2):249-274, 2018.

\bibitem{Razborov}
A. Razborov, Flag algebras.
{\em Journal of Symbolic Logic}, 72:1239--1282, 2007.

\bibitem{Riener}
C. Riener.
On the degree and half degree principle for symmetric polynomials.
{\em Journal of Pure and Applied Algebra}, 216(4):850--856, 2012.

\bibitem{RTAL}
C. Riener, T. Theobald, L. Jansson Andr\'en, and J.B. Lasserre.
Exploiting symmetry in SDP-relaxations for polynomial optimization.
{\em Mathematics of Operations Research}, 38(1):122--141, 2013.

\bibitem{Schrijver1979}
A. Schrijver.
A comparison of the Delsarte and Lov\'asz bounds.
{\em IEEE Transactions on Information Theory}, 25:425--429, 1979.

\bibitem{Schrijver2005}
A. Schrijver.
 New code upper bounds from the Terwilliger algebra and semidefinite programming.
  {\em IEEE Transactions on Information Theory}, 51:2859--2866, 2005.










\bibitem{ArtinWedderburn} 
J.H.M.  Wedderburn.  {\em Lectures on matrices.}
Dover Publications Inc., New York, 1964.
\end{thebibliography}
\end{document}